%% file: TLSSfinal.tex
\newif\ifprint%IM
\printfalse%IM off for contributors
\documentclass[twoside,10pt]{HandbookOfModuli}%IM

\usepackage{amsmath,amsthm}
\usepackage{amssymb}
\usepackage{latexsym}%IM
\usepackage[utf8]{inputenc}%IM
\usepackage{textcomp}%IM
\usepackage{color}%IM
\usepackage{titlesec}
\usepackage{xspace}

\ifprint
	\input{giovanni} % off for contributors
	\usepackage[final]{microtype}%IM off for contributors
\fi
\usepackage{bm}
\renewcommand{\mathbf}[1]{\bm{#1}} % override \mathbf for use with Giovanni

\ifprint
	\definecolor{linkred}{rgb}{0,0,0} % black
	\definecolor{linkblue}{rgb}{0,0,0} % black
\else
	\definecolor{linkred}{rgb}{0.7,0.2,0.2}
	\definecolor{linkblue}{rgb}{0,0.2,0.6}
\fi

\numberwithin{equation}{section} % standardize numbering

\catcode`@=11 \baselineskip 4.5mm
\def\ps@handbook{\def\@oddhead{\hfill \leftmark \hfill\thepage }
\def\@evenhead{\thepage \hfill \rightmark \hfill}
\def\@oddfoot{}
\def\@evenfoot{}}
\def\@evenhead{}
\def\@oddfoot{}%\hfill \copyright\ China Higher Education Press}
\def\@evenfoot{\hfill\copyright\ China Higher Education Press}
\def\list#1#2{\ifnum \@listdepth >5\relax \@toodeep \else \global
\advance \@listdepth\@ne \fi \rightmargin \z@ \listparindent\z@
\itemindent\z@ \csname @list\romannumeral\the\@listdepth\endcsname
\def\@itemlabel{#1}\let\makelabel\@mklab \@nmbrlistfalse #2\relax
\@trivlist \parskip -\parsep \parindent\listparindent \advance
\linewidth -\rightmargin \advance\linewidth -\leftmargin \advance
\@totalleftmargin \leftmargin \parshape \@ne \@totalleftmargin
\linewidth \ignorespaces}
\renewcommand*\l@section{\@tocline{1}{0pt}{0em}{1em}{}}%IM 1/10 fix toc incompatibility
\renewcommand*\l@subsection{\@tocline{2}{0pt}{1.5em}{2em}{}} %IM 1/10 fix toc incompartibility 
\renewcommand\contentsnamefont{\bfseries\large} %IM 1/10 fix toc incompartibility 
\catcode`@=12
\renewcommand{\theequation}{\thesection.\arabic{equation}}

\def\thebibliography#1{\section*{References}
\list{[\arabic{enumi}]}{\settowidth \labelwidth{[#1]} \leftmargin
\labelwidth \advance \leftmargin \labelsep \usecounter{enumi}}
\def\newblock{\hskip .11em plus .33em minus .07em} \sloppy
\clubpenalty 4000 \widowpenalty 4000 \sfcode`\.=1000 \relax}
\titleformat{\section}{\normalfont\large\bfseries}{\thesection.}{0.5em}{}[\kern0.em]%IM 3/3 null kern marks mark optional vertical spacing after
\titleformat{\subsection}{\normalfont\bfseries}{\thesubsection.}{0.3em}{}[\kern0.em]%IM 3/3 null kern marks mark optional vertical spacing after
\titleformat{\subsubsection}[runin]{\normalfont\bfseries}{\thesubsubsection.}{0.5em}{}[\kern0.5em]%IM 9/29 set as run in head

\def\fofsubsubsection#1{\refstepcounter{equation}\subsubsection*{\theequation.\kern0.25em #1}}%IM 3/3 flush "subsubsection" using section.equation numbering
\def\foisubsubsection#1{\refstepcounter{equation}\subsubsection*{\kern\parindent\theequation.\kern0.25em #1}}%IM 3/3 indented "subsubsection" using section.equation numbering   

\usepackage{amscd}
\usepackage{array}

\input{TorelliMacros}
\usepackage{url}

\newtheorem{theorem}[equation]{Theorem}
\newtheorem{proposition}[equation]{Proposition}
\newtheorem{corollary}[equation]{Corollary}
\newtheorem{conjecture}[equation]{Conjecture}
\newtheorem{expectation}[equation]{Expectation}

\theoremstyle{definition}
\newtheorem{definition}[equation]{Definition}
\newtheorem{example}[equation]{Example}
\newtheorem{remark}[equation]{Remark}

\newtheorem{question}[equation]{Question}

\numberwithin{equation}{section}
\newcommand{\mosubsection}[1]{\setcounter{subsection}{\value{equation}}\stepcounter{equation}\subsection{#1}}

\begin{document}
\setcounter{page}{1}% to be adjusted when all contributions are assembled

%%%%%%%%%%%%%%%%%%%%%%%%%%%%%%%%%%%%%%%%%%%%%%%%%%%%%%%%%%%%
%%
%% Start of the title page
%%
%%%%%%%%%%%%%%%%%%%%%%%%%%%%%%%%%%%%%%%%%%%%%%%%%%%%%%%%%%%%

\title{The Torelli locus and special subvarieties}
   
\author{Ben Moonen}
\address{University of Amsterdam}
\email{bmoonen@uva.nl}

\author{Frans Oort}
\address{University of Utrecht}
\email{f.oort@uu.nl}
   

\begin{abstract}
We study the Torelli locus $\Torelli_g$ in the moduli space~$\sfA_g$ of abelian varieties. We consider special subvarieties (Shimura subvarieties) contained in the Torelli locus. We review the construction of some non-trivial examples, and we discuss some conjectures, techniques and recent progress.
\end{abstract}

\maketitle
\thispagestyle{empty}% suppress page number on title pagesr

%%%%%%%%%%%%%%%%%%%%%%%%%%%%%%%%%%%%%%%%%%%%%%%%%%%%%%%%%%%%
%%
%% Start of the body of the article
%%
%%%%%%%%%%%%%%%%%%%%%%%%%%%%%%%%%%%%%%%%%%%%%%%%%%%%%%%%%%%%

\section*{Introduction}

\subsection*{The Torelli morphism}
An algebraic curve~$C$ has a Jacobian $J_C = \Pic^0(C)$, which comes equipped with a canonical principal polarization $\lambda \colon J_C \isomarrow J_C^t$. The dimension of~$J_C$ equals the genus~$g$ of the curve. One should like to understand~$J_C$ knowing certain properties of~$C$, and, conversely, we like to obtain information about~$C$ from properties of~$J_C$ or the pair $(J_C,\lambda)$.

As an important tool for this we have the moduli spaces of
algebraic curves and of polarized abelian varieties. The construction that associates to~$C$ the principally polarized abelian variety $(J_C,\lambda)$ defines a morphism
$$
j\colon \sfM_g \to \sfA_g\, ,
$$
the \emph{Torelli morphism}, which by a famous theorem of Torelli is injective on geometric points. (Here we work with the coarse moduli schemes.) The image  $\Torelli^\open_g := j(\sfM_g)$ is called the (open) Torelli locus, also called the Jacobi locus. Its Zariski closure $\Torelli_g$ inside~$\sfA_g$ is called the \emph{Torelli locus}.  We should like to study the interplay between curves and their Jacobians using the geometry of $\sfM_g$ and of~$\sfA_g$, and relating these via the inclusion $\Torelli_g \subset \sfA_g$.

\subsection*{Special subvarieties}
Although the moduli spaces $\sfM_g$ and~$\sfA_g$ have been studied extensively, it turns out that there are many basic questions that are still open. The moduli space~$\sfA_g$ locally has a group-like structure, in a sense that can be made precise. As a complex manifold, $\sfA_g(\bC)$ is a quotient of the Siegel space, which is a principally homogeneous space under the group $\CSp_{2g}(\bR)$. Hence it is natural to study subvarieties that appear as images of orbits of an algebraic subgroup;
such a subvariety is called a \emph{special subvariety}. (We refer to Section~\ref{SpecSubvars} for a more precise definition of this notion.)

The zero dimensional special subvarieties are precisely the CM points, corresponding to abelian varieties where $\End(A) \otimes \bQ$ contains a commutative semi-simple algebra of rank $2{\cdot}\dim(A)$. These are fairly well understood. The arithmetical properties of CM points form a rich and beautiful subject, that plays an important role in the study of the moduli space. The presence of a dense set of CM points enabled Shimura to prove for a large class of Shimura varieties that they admit a model over a number field. In its modern form this theory owes much to Deligne's presentation in \cite{Del-TravShim} and~\cite{Del-VarSh}, and the existence of canonical models for all Shimura varieties was established by Borovoi and Milne.

For the problems we want to discuss, the special points again play a key role. It is not hard to show that on any special subvariety of~$\sfA_g$ over~$\bC$ the special points are dense, even for the analytic topology. A conjecture by Y.~Andr\'e and F.~Oort says that, conversely, an algebraic subvariety with a
Zariski dense set of CM points is in fact a special variety. See \ref{AOConj}. There is a nice analogy between this conjecture and the Manin-Mumford conjecture, proven by Raynaud, which says that an algebraic subvariety of an abelian variety that contains a Zariski dense set of torsion points is the translate of an abelian subvariety under a torsion point.

Klingler and Yafaev, using work of Ullmo and Yafaev, have announced a proof of the Andr\'e-Oort conjecture under assumption of the Generalized Riemann Hypothesis for CM fields.

\subsection*{Coleman's Conjecture and special subvarieties in the Torelli locus} 
In view of the above, it is natural to investigate the subvariety $\Torelli_g \subset \sfA_g$ from the perspective of special subvarieties of~$\sfA_g$.
A conjecture by Coleman, again based on the analogy with the Manin-Mumford conjecture, says that for $g>3$ there are only finitely many complex curves of genus~$g$ such that $J_C$ is an abelian variety with CM. See~\ref{ColemanConj}. In terms of moduli spaces, this says that for $g>3$ the number of special points on $\Torelli^\open_g$ should be finite. Combining this with the Andr\'e-Oort conjecture, it suggests that for large~$g$ there are no special subvarieties $Z \subset \sfA_g$ of positive dimension that are contained inside~$\Torelli_g$ and meet~$\Torelli^\open_g$; see~\ref{ZinTorelli}. We should like to point out that there is no clear evidence for this expectation, although there are several partial results in support of it. See for instance Theorem~\ref{WeylCM} below, which gives a weaker version of Coleman's conjecture.

As we shall explain in Section~\ref{ExaSsT}, Coleman's Conjecture is not true for $g \in \{4,5,6,7\}$. The reason is that for these genera we can find explicit families of curves such that the corresponding Jacobians trace out a special subvariety of positive dimension in~$\sfA_g$. All known examples of this kind, for $g \geq 4$, arise from families of abelian covers of~$\bP^1$, where we fix the Galois group and the local monodromies, and we vary the branch points. For cyclic covers it can be shown that, beyond the known examples, this construction gives no further families that give rise to a special subvariety in~$\sfA_g$; see Theorem~\ref{BM20Fams}.

In this paper, we give a survey of what is known about special subvarieties in the Torelli locus, and we try to give the reader a feeling for the difficulties one encounters in studying these. Going further, one could ask for a complete classification of such special subvarieties $Z \subset \Torelli_g$ for every~$g$. We should like to point out that already for $g=4$ we do not have such a classification, in any sense.

\subsection*{Acknowledgments}
We thank Dick Hain for patiently explaining to us some details related to his paper~\cite{Hain-Locsymm}.

\mosubsection{Notation}
\label{Notation}

If $G$ is a group scheme of finite type over a field~$k$, we denote its identity component by~$G^0$. In case $k$ is a subfield of~$\bR$ we write $G(\bR)^+$ for the identity component of $G(\bR)$ with regard to the analytic topology. If $G$ is reductive, $G^\ad$ denotes the adjoint group.

\emph{The group of symplectic similitudes.} Let a positive integer~$g$ be given. Write $V := \bZ^{2g}$, and let $\phi \colon V \times V \to \bZ$ be the standard symplectic pairing, so that on the standard ordered basis~$\{e_i\}_{i=1,\ldots,2g}$ we have $\phi(e_i,e_j) = 0$ if $i+j \neq 2g+1$ and $\phi(e_i,e_{2g+1-i}) = 1$ if $i \leq g$. The group $\CSp(V,\phi)$ is the reductive group over~$\bZ$ of symplectic similitudes of~$V$ with regard to the form~$\phi$. If there is no risk of confusion we simply write $\CSp_{2g}$ for this group. We denote by $\nu \colon \CSp_{2g} \to \bG_\mult$ the multiplier character. The kernel of~$\nu$ is the symplectic group~$\Sp_{2g}$ of transformations of~$V$ that preserve the form~$\phi$.

\setcounter{tocdepth}{1}
\tableofcontents

\section{The Torelli morphism}

\mosubsection{Moduli of abelian varieties}
Fix an integer $g > 0$. If $S$ is a base scheme and $A$ is an abelian scheme over~$S$, we write~$A^t$ for the dual abelian scheme. There is a canonical homomorphism $\kappa_A\colon A \to A^{tt}$, which Cartier and Nishi proved to be an isomorphism; see \cite{Cartier1}, \cite{Cartier2}, \cite{Nishi}; also see \cite{FO-CGS}, Theorem~20.2. 

On $A \times_S A^t$ we have a Poincar\'e bundle~$\cP$. A \emph{polarization} of~$A$ is a symmetric isogeny $\lambda \colon A \to A^t$ with the property that the pull-back of~$\cP$ via the morphism $(\id,\lambda) \colon A \to A \times_S A^t$ is a relatively ample bundle on~$A$ over~$S$. The polarization~$\lambda$ is said to be principal if in addition it is an isomorphism of abelian schemes over~$S$. A pair $(A,\lambda)$ consisting of an abelian scheme over~$S$ together with a principal polarization is called a \emph{principally polarized abelian scheme}. We abbreviate ``principally polarized abelian scheme(s)/variety/varieties'' to ``ppav'' (sic!). 

We denote by~$\mathcal{A}_g$ the moduli stack over~$\Spec(\mathbb{Z})$ of $g$-dimensional ppav. It is a connected Deligne-Mumford stack that is quasi-projective and smooth over~$\mathbb{Z}$, of relative dimension $g(g+1)/2$.  The characteristic zero fiber $\mathcal{A}_{g,\bQ}$ can be shown to be geometrically irreducible by transcendental methods; cf.\ Section~\ref{AgShimVar}. Using this, Chai and Faltings have proved that the characteristic~$p$ fibers are geometrically irreducible, too; see \cite{Faltings-Bonn}, \cite{ChaiPhD}, and \cite{F.C} Chap.~IV, Corollary~6.8. A pure characteristic~$p$ proof can be obtained using the results of~\cite{FO-EO}.

Note that later in this article we shall use the same notation~$\mathcal{A}_g$ for the moduli stack over some given base field.

For many purposes it is relevant to consider ppav together with a level structure. We introduce some notation and review some facts. Let again $(A,\lambda)$ be a ppav of relative dimension~$g$ over some base scheme~$S$. Given $m \in \mathbb{Z}$ we have a ``multiplication by~$m$'' endomorphism $[m] \colon A\to A$. We denote by $A[m]$ the kernel of this endomorphism, in the scheme-theoretic sense. Assume that $m \neq 0$. Then $A[m]$ is a finite locally free group scheme of rank~$m^{2g}$. There is a naturally defined bilinear pairing
\[
e_m^\lambda \colon A[m] \times A[m] \to \mu_m\, ,
\]
called the Weil pairing, that is non-degenerate and that is symplectic in the sense that $e_m^\lambda(x,x) = 1$ for all $S$-schemes~$T$ and all $T$-valued points $x \in A[m]\bigl(T\bigr)$. If $m$ is invertible in $\Gamma(S,O_S)$ then the group scheme $A[m]$ is \'etale over~$S$.

Fix an integer $m \geq 1$. Write $\bZ[\zeta_m]$ for $\bZ[t]/(\Phi_m)$, where $\Phi_m \in \bZ[t]$ is the $m$th cyclotomic polynomial. We shall first define a Deligne-Mumford stack $\mathcal{A}^\prime_{g,[m]}$ over the spectrum of the ring $R_m := \bZ[\zeta_m, 1/m]$. The moduli stack $\mathcal{A}_{g,[m]}$ is then defined to be the same stack, but now viewed as a stack over $\bZ[1/m]$, with as structural morphism the composition 
\[
\mathcal{A}^\prime_{g,[m]} \to \Spec(R_m) \to \Spec\bigl(\bZ[1/m]\bigr)\, .
\]
The reason that we want to view $\mathcal{A}_{g,[m]}$ as an algebraic stack over $\bZ[1/m]$ is that this is the stack whose characteristic zero fiber has a natural interpretation as a Shimura variety; see the discussion in~\ref{AgShimVar} below. 

Let $S$ be a scheme over~$R_m$. Note that this just means that $S$ is a scheme such that $m$ is invertible in $\Gamma(S,O_S)$, and that we are given a primitive $m$th root of unity $\zeta \in \Gamma(S,O_S)$. We may view this~$\zeta$ as an isomorphism of group schemes $\zeta \colon (\bZ/m\bZ) \isomarrow \mu_m$. Let $V$ and~$\phi$ be as defined in~\ref{Notation}. Given a ppav $(A,\lambda)$ of relative dimension~$g$ over~$S$, by a \emph{(symplectic) level~$m$ structure} on $(A,\lambda)$ we then mean an isomorphism of group schemes $\alpha \colon (V/mV) \isomarrow A[m]$ such that
the diagram
\begin{equation*}
\begin{CD}
(V/mV) \times (V/mV) @>\phi>> (\bZ/m\bZ)\\
@V{\alpha\times\alpha}VV @VV{\zeta}V\\
A[m] \times A[m] @>e_m^\lambda>> \mu_m
\end{CD}
\end{equation*}
is commutative. 

We write $\mathcal{A}^\prime_{g,[m]}$ for the moduli stack over~$R_m$ of ppav with a level~$m$ structure. Again this is a quasi-projective smooth Deligne-Mumford stack with geometrically irreducible fibers. If $m \geq 3$, it is a scheme. The resulting stack
\[
\cA_{g,[m]} := \Bigl( \cA^\prime_{g,[m]} \to \Spec(R_m) \to \Spec\bigl(\bZ[1/m]\bigr)\Bigr)
\]
is a quasi-projective smooth Deligne-Mumford stack over~$\bZ[1/m]$ whose geometric fibers have $\varphi(m)$ irreducible components.

We now assume that $m \geq 3$. The finite group $\mathcal{G} := \CSp_{2g}(\bZ/m\bZ)$ acts on the stack $\cA^\prime_{g,[m]}$ via its action on the level structures. Note that this is not an action over~$R_m$; only the subgroup $\Sp_{2g}(\bZ/m\bZ)$ acts by automorphisms over~$R_m$. We do, however, obtain an induced action of~$\mathcal{G}$ on~$\cA_{g,[m]}$ over~$\bZ[1/m]$. The stack quotient of $\cA_{g,[m]}$ modulo~$\mathcal{G}$ is just $\cA_g \otimes \bZ[1/m]$. The scheme quotient of $\cA_{g,[m]}$ modulo~$\mathcal{G}$ is the coarse moduli scheme $\sfA_g \otimes \bZ[1/m]$. For $m_1$, $m_2 \geq 3$ these coarse moduli schemes agree over $\bZ[1/(m_1m_2)]$; hence we can glue them to obtain a coarse moduli scheme~$\sfA_g$ over~$\bZ$.

The stacks $\cA_{g,[m]}$ are not complete. There are various ways to compactify them. The basic reference for this is the book~\cite{F.C}. There are toroidal compactifications that depend on the choice of some cone decomposition~$\Sigma$ (see op.\ cit.\ for the details); these compactifications map down to the Baily-Borel (or minimal, or Satake) compactification.

\mosubsection{Moduli of curves}
Fix an integer $g \geq 2$. We denote by~$\cM_g$ the moduli stack of curves of genus~$g$ over~$\Spec(\bZ)$. It is a quasi-projective smooth Deligne-Mumford stack over~$\bZ$ of relative dimension $3g-3$ for $g\geq 2$, with geometrically irreducible fibers.

We write $\ol\cM_g$ for the Deligne-Mumford compactification of~$\cM_g$; it is the stack of stable curves of genus~$g$. See~\cite{Del.Mumf}. The boundary $\ol\cM_g\setminus \cM_g$ is a divisor with normal crossings; it has irreducible components~$\Delta_i$ for $0 \leq i \leq \lfloor g/2\rfloor$. There is a non-empty open part of~$\Delta_0$ that parametrizes irreducible curves with a single node. Similarly, for $i>0$ there is a non-empty open part of~$\Delta_i$ parametrizing curves with precisely two irreducible components, of genera~$i$ and~$g-i$. The complement of~$\Delta_0$ is the open substack $\cM_g^\ct \subset \ol\cM_g$ of curves of compact type. Here we recall that a stable curve~$C$ over a field~$k$ is said to be \emph{of compact type} if the connected components of its Picard scheme are proper over~$k$. If $k$ is algebraically closed this is equivalent to the condition that the dual graph of the curve has trivial homology; so the irreducible components of~$C$ are regular curves, the graph of components is a tree, and the sum of the genera of the irreducible components equals~$g$. Note that $\cM_g \subset \cM_g^\ct$.

For the corresponding coarse moduli schemes we use the notation $\sfM_g$ and~$\sfM_g^\ct$.

\mosubsection{The Torelli morphism and the Torelli locus}
\label{TorelliLocus}
Let $C$ be a stable curve of compact type over a base scheme~$S$, with fibers of genus $g \geq 2$. The Picard scheme $\Pic_{C/S}$ is a smooth separated $S$-group scheme whose components are of finite type and proper over~$S$. In particular, the identity component $J_C := \Pic^0_{C/S}$ is an abelian scheme over~$S$. It has relative dimension~$g$ and comes equipped with a canonical principal polarization~$\lambda$.

The functor that sends $C$ to the pair $(J_C,\lambda)$ defines a representable morphism of algebraic stacks
\[
j \colon \cM_g^\ct \to \cA_g\, ,
\]
called the \emph{Torelli morphism}. It is known that this morphism is proper. 

Torelli's celebrated theorem says that for an algebraically closed field~$k$ the restricted morphism $j \colon \cM_g \to \cA_g$ is injective on $k$-valued points. In other words, if $C_1$ and~$C_2$ are smooth curves over~$k$ such that $\bigl(J_{C_1},\lambda_1\bigr)$ and $\bigl(J_{C_2},\lambda_2\bigr)$ are isomorphic then $C_1$ and~$C_2$ are isomorphic.  This injectivity property of the Torelli morphism does not hold over the boundary; for instance, if $C$ is a stable curve of compact type with two irreducible components then the Jacobian of~$C$ does not ``see'' at which points the two components are glued. For a detailed study of what happens over the boundary, see~\cite{Caporaso.Viviani}.

Let us now work over a field~$k$. On coarse moduli schemes the Torelli morphism gives rise to a morphism $j \colon \sfM_g^\ct \to \sfA_g$. We define the \emph{Torelli locus} $\Torelli_g \subset \sfA_g$ to be the image of this morphism. It is a reduced closed subscheme of~$\sfA_g$. By the \emph{open Torelli locus} $\Torelli_g^\open \subset \Torelli_g$ we mean the image of~$\sfM_g$, i.e., the subscheme of~$\Torelli_g$ whose geometric points correspond to Jacobians of nonsingular curves. Note that $\Torelli_g^\open \subset \Torelli_g$ is open and dense. 

The boundary $\Torelli_g\setminus \Torelli_g^\open$ is denoted by~$\Torelli_g^\dec$. The notation refers to the fact that a point $s \in \Torelli_g(\overline{k})$ is in~$\Torelli_g^\dec$ if and only if the corresponding ppav $(A_s,\lambda_s)$ is decomposable, meaning that we have ppav $(B_1,\mu_1)$ and $(B_2,\lambda_2)$ of smaller dimension such that $(A_s,\lambda_s) \cong (B_1,\mu_1) \times (B_2,\mu_2)$ as ppav. In one direction this is clear, for if $C$ is a reducible curve of compact type, $(J_C,\lambda)$ is the product, as a ppav, of the Jacobians of the irreducible components of~$C$. Conversely, if $(A_s,\lambda_s)$ is decomposable, it has a symmetric theta divisor that is reducible, which implies it cannot be the Jacobian of an irreducible (smooth and proper) curve.

\begin{remark}
Let $g>2$. The Torelli morphism $j \colon \cM_g \to \cA_g$ is ramified at the hyperelliptic locus. Outside the hyperelliptic locus it is an immersion. The picture is different for the Torelli morphism $j \colon \sfM_g \to \sfA_g$ on coarse moduli schemes. The morphism $j \colon \sfM_{g,\bQ} \to \sfA_{g,\bQ}$ on the characteristic zero fibers is an immersion; however, in positive characteristic this is not true in general. See \cite{FO-JHMS}, Corollaries 2.8, 3.2 and~5.3.
\end{remark}

\section{Some abstract Hodge theory}

In this section we review some basic notions from abstract Hodge theory. We shall focus on examples related to abelian varieties.

\mosubsection{Basic definitions}
\label{HodgeTh}
We start by recalling some basic definitions. For a more comprehensive treatment we refer to \cite{Del-HodgeII} and~\cite{Peters-Steenbr}.

Let $\bS := \Res_{\bC/\bR}(\Gm)$, which is an algebraic torus over~$\bR$ of rank~$2$, called the Deligne torus. Its character group is isomorphic to~$\bZ^2$, on which complex conjugation, the non-trivial element of $\Gal(\bC/\bR)$, acts by $(m_1,m_2) \mapsto (m_2,m_1)$. We have $\bS(\bR) = \bR^*$ and $\bS(\bC) = \bC^* \times \bC^*$. Let $w \colon \bG_{\mult,\bR} \to \bS$ be the cocharacter given on real points by the natural embedding $\bR^* \hookrightarrow \bC^*$; it is called the weight cocharacter. Let $\Nm\colon \bS \to \bG_{\mult,\bR}$ be the homomorphism given on real points by $z\mapsto z\bar{z}$; it is called the norm character.

A \emph{$\bZ$-Hodge structure of weight~$k$} is a torsion-free $\bZ$-module~$H$ of finite rank, together with a homomorphism of real algebraic groups $h \colon \bS \to \GL(H)_\bR$ such that $h \circ w \colon \Gm \to \GL(H)_\bR$ is the homomorphism given by $z \mapsto z^{-k} \cdot \id_H$. The Hodge decomposition
\[
H_\bC = \oplus_{p+q=k}\, H_\bC^{p,q}
\]
is obtained by taking $H_\bC^{p,q}$ to be the subspace of~$H_\bC$ on which $(z_1,z_2) \in \bC^* \times \bC^* = \bS(\bC)$ acts as multiplication by $z_1^{-p} z_2^{-q}$. (We follow the sign convention of~\cite{Del-HConAV}; see loc.\ cit., Remark~3.3.) 

As an example, the Tate structure $\bZ(n)$ is the Hodge structure of weight $-2n$ with underlying $\bZ$-module $\bZ(n) = (2\pi i)^n \cdot \bZ$ given by the homomorphism $\Nm^n \colon \bS \to \bG_{\mult,\bR} = \GL\bigl(\bZ(n)\bigr)_\bR$. If $H$ is any Hodge structure, we write $H(n)$ for $H \otimes_\bZ \bZ(n)$.

The endomorphism $C = h(i)$ of~$H_\bR$ is known as the Weil operator. If $H$ is a $\bZ$-Hodge structure of weight~$k$, a \emph{polarization} of~$H$ is a homomorphism of Hodge structures $\psi \colon H \otimes H \to \bZ(-k)$ such that the bilinear form
\[
H_\bR \times H_\bR \to \bR
\]
given by $(x,y) \mapsto (2\pi i)^k \cdot \psi_\bR\bigl(x \otimes C(y)\bigr)$ is symmetric and positive definite. (Here $\psi_\bR$ is obtained from~$\psi$ by extension of scalars to~$\bR$.) This condition implies that $\psi$~is symmetric if $k$ is even and alternating if $k$ is odd.

Instead of working with integral coefficients, we may also consider $\bQ$-Hodge structures. The above definitions carry over verbatim.

\begin{example}
Let $(A,\lambda)$ be a polarized complex abelian variety. The first homology group $H = H_1(A,\bZ)$ carries a canonical Hodge structure of type $(-1,0) + (0,-1)$. (By this we mean that $H^{p,q}= (0)$ for all pairs $(p,q)$ different from $(-1,0)$ and $(0,-1)$.) The polarization $\lambda \colon A \to A^t$ (in the sense of the theory of abelian varieties) induces a polarization (in the sense of Hodge theory) $\psi \colon H \otimes H \to \bZ(1)$. Let $\CSp(H,\psi) \subset \GL(H)$ be the algebraic group of symplectic similitudes of~$H$ with respect to the symplectic form~$\psi$. The homomorphism $h \colon \bS \to \GL(H)_\bR$ that gives the Hodge structure factors through $\CSp(H,\psi)_\bR$.

A crucial fact for much of the theory we want to discuss is that the map $(A,\lambda) \mapsto (H,\psi)$ gives an equivalence of categories
\begin{equation}
\label{AVEqCat}
\left\{\vcenter{
\setbox0=\hbox{polarized complex}
\copy0
\hbox to\wd0{\hfill  abelian varieties\hfill}}\right\}
\xrightarrow{\; \eq\;}
\left\{\vcenter{
\setbox0=\hbox{polarized Hodge structures}
\copy0
\hbox to\wd0{\hfill of type $(-1,0) + (0,-1)$\hfill}}
\right\}\; .
\end{equation}
In \ref{AgShimVar} below we shall explain how this fact gives rise to a modular interpretation of certain Shimura varieties. As a variant, the category of abelian varieties is equivalent to the category of polarizable $\bZ$-Hodge structures. Another variant is obtained by working with $\bQ$-coefficients; we get that the category of abelian varieties up to isogeny is equivalent to the category of polarizable $\bQ$-Hodge structures of type $(-1,0) + (0,-1)$.
\end{example}

\mosubsection{Hodge classes and Mumford-Tate groups}
\label{MTGps}
Let $H$ be a $\bQ$-Hodge structure of weight~$k$. By a \emph{Hodge class} in~$H$ we mean an element of
\[
H \cap H_\bC^{0,0}\, ;
\]
so, a \emph{rational} class that is purely of type $(0,0)$ in the Hodge decomposition. Clearly, a non-zero Hodge class can exist only if $k=0$. If the weight is~$0$ and $\Fil^\bullet H_\bC$ is the Hodge filtration, the space of Hodge classes is also equal to
\[
H \cap \Fil^0 H_\bC\, ,
\]
where again the intersection is taken inside~$H_\bC$.

As an example, we note that a $\bQ$-linear map $H \to H$, viewed as an element of $\End_\bQ(H) = H^\vee \otimes H$, is a homomorphism of Hodge structures if and only it is a Hodge class for the induced Hodge structure on $\End_\bQ(H)$. (Note that $H^\vee \otimes H$ indeed has weight~$0$.) Similarly, a polarization $H \otimes H \to \bQ(-k)$ gives rise to a Hodge class in $(H^\vee)^{\otimes 2}(-k)$. 

Let $h \colon \bS \to \GL(H)_\bR$ be the homomorphism that defines the Hodge structure on~$H$. Consider the homomorphism $(h,\Nm) \colon \bS \to \GL(H)_\bR \times \bG_{\mult,\bR}$. The \emph{Mumford-Tate group of~$H$}, notation $\MT(H)$, is the smallest algebraic subgroup $M \subset \GL(H) \times \Gm$ over~$\bQ$ such that $(h,\Nm)$ factors through~$M_\bR$. As we shall see in Example~\ref{MTofAV} below, in some important cases the definition takes a somewhat simpler form.

It is known that $\MT(H)$ is a connected algebraic group. If the Hodge structure~$H$ is polarizable, $\MT(H)$ is a reductive group; this plays an important role in many applications.

The crucial property of the Mumford-Tate group is that it allows us to calculate spaces of Hodge classes, as we shall now explain. For simplicity we shall assume the Hodge structure~$H$ to be polarizable.

We consider the Hodge structures we can build from~$H$ by taking direct sums, duals, tensor products and Tate twists. Concretely, given a triple $\mathbf{m} = (m_1,m_2,m_3)$ with $m_1$, $m_2 \in \bZ_{\geq 0}$ and $m_3 \in \bZ$, define
\[
T(\mathbf{m}) := H^{\otimes m_1} \otimes (H^\vee)^{\otimes m_2} \otimes \bQ(m_3)\, .
\]
The group $\GL(H) \times \Gm$ naturally acts on~$H$ via the projection onto its first factor; hence it also acts on $H^\vee$. We let it act on $\bQ(1)$ through the second projection. This gives us an induced action of $\GL(H) \times \Gm$ on~$T(\mathbf{m})$. The crucial property of the Mumford-Tate group $\MT(H) \subset \GL(H) \times \Gm$ is that for any such tensor construction $T = T(\mathbf{m})$, the subspace of Hodge classes $T \cap T_\bC^{0,0}$ is precisely the subspace of $\MT(H)$-invariants. Note that the $\MT(H)$-invariants are just the invariants in~$T$ under the action of the group $\MT(H)\bigl(\bQ\bigr)$ of $\bQ$-rational points of the Mumford-Tate group. As the $\bQ$-points of $\MT(H)$ are Zariski dense in the $\bC$-group $\MT(H)_\bC$, the subspace of $\MT(H)\bigl(\bC\bigr)$-invariants in $T \otimes \bC$ equals
\[
\bigl(T \cap T_\bC^{0,0}\bigr) \otimes_\bQ \bC\, .
\]

The above property of the Mumford-Tate group characterizes it uniquely. In other words, if $M \subset \GL(H) \times \Gm$ is an algebraic subgroup with the property that for any $T = T(\mathbf{m})$ the space of Hodge classes $T \cap T_\bC^{0,0}$ equals the subspace of $M$-invariants in~$T$, we have $M = \MT(H)$. (See~\cite{Del-HConAV}, Section~3.) In this way we have a direct coupling between the Mumford-Tate group and the various spaces of Hodge classes in tensor constructions obtained from~$H$.

\begin{example}
\label{MTofAV}
Consider a polarized abelian variety $(A,\lambda)$ over~$\bC$. Let $H_\bQ=H_1(A,\bQ)$ with polarization form $\psi \colon H_\bQ \otimes H_\bQ \to \bQ(1)$ be the associated $\bQ$-Hodge structure, and recall that the Hodge structure is given by a homomorphism $h \colon \bS \to \CSp(H_\bQ,\psi)_\bR$. Define the Mumford-Tate group of~$A$, notation $\MT(A)$, to be the smallest algebraic subgroup $M \subset \CSp(H_\bQ,\psi)$ over~$\bQ$ with the property that $h$ factors through~$M_\bR$.

The Mumford-Tate group as defined here is really the same as the Mumford-Tate group of the Hodge structure~$H_\bQ$. To be precise, $\MT(H_\bQ) \subset \GL(H_\bQ) \times \Gm$ as defined above is the graph of the multiplier character $\nu \colon \CSp(H_\bQ,\psi) \to \Gm$ restricted to~$\MT(A)$. When considering Hodge classes in tensor constructions~$T$ as above, we let $\MT(A)$ act on the Tate structure~$\bQ(1)$ through the multiplier character.

If we want to understand how the Mumford-Tate group is used, the simplest non-trivial examples are obtained by looking at endomorphisms. As customary we write $\End^0(A)$ for $\End(A) \otimes \bQ$. By the equivalence of categories of~(\ref{AVEqCat}), in the version with $\bQ$-coefficients, we have
\begin{equation}
\label{End=End}
\End^0(A) \xrightarrow{\;\sim\;} \End_{\QHS}(H_\bQ) = \End_\bQ(H_\bQ)^{\MT(A)}\, .
\end{equation}
Once we know the Mumford-Tate group, this allows us to calculate the endomorphism algebra of~$A$. In practice we often use this the other way around, in that we assume $\End^0(A)$ known and we use~(\ref{End=End}) to obtain information about $\MT(A)$.

In general Hodge classes on~$A$ are not so easy to interpret geo\-metrically; see for instance the example in \cite{Mumf-ShimPaper}, \S~4. (In this example, $A$ is a simple abelian fourfold and there are Hodge classes in $H^4(A^2,\bQ)$ that are not in the algebra generated by divisor classes; see \cite{BM.YuZ-Duke} and \cite{BM.YuZ-MathAnn}, (2.5).) According to the Hodge conjecture, all Hodge classes should arise from algebraic cycles on~$A$ but even for abelian varieties this is far from being proven.
\end{example}

\mosubsection{Hodge loci}
\label{HodgeLoci}
We study the behavior of Mumford-Tate groups in families. The setup here is that we consider a polarizable $\bQ$-VHS (variation of Hodge structure) of weight~$k$ over a complex manifold~$S$. Let $\mathcal{H}$ denote the underlying $\bQ$-local system. For each $s \in S(\bC)$ we have a Mumford-Tate group $\MT_s \subset \GL\bigl(\mathcal{H}(s)\bigr) \times \Gm$, and we should like to understand how this group varies with~$s$. This is best described by passing to a universal cover $\pi \colon \tilde{S} \to S$. Write $\tilde{\mathcal{H}} := \pi^* \mathcal{H}$, and let $H := \Gamma(\tilde{S},\tilde{\mathcal{H}})$. We have a trivialization $\tilde{\mathcal{H}} \isomarrow H \times \tilde{S}$. Using this we may view the Mumford-Tate group $\MT_x$ of a point $x \in \tilde{S}$ (by which we really mean the Mumford-Tate group of the image of~$x$ in~$S$) as an algebraic subgroup of $\GL(H) \times \Gm$. So, passing to the universal cover has the advantage that we can describe the VHS as a varying family of Hodge structures on some given space~$H$ and that we may describe the Mumford-Tate groups $\MT_x$ as subgroups of one and the same group.

Given a triple $\mathbf{m} = (m_1,m_2,m_3)$ as in~\ref{MTGps}, consider the space
\[
T(\mathbf{m}) := H^{\otimes m_1} \otimes (H^\vee)^{\otimes m_2} \otimes \bQ(m_3)\, ,
\]
on which we have a natural action of $\GL(H) \times \Gm$. We may view an element $t \in T(\mathbf{m})$ as a global section of the local system
\[
\tilde{\mathcal{T}}(\mathbf{m}) := \tilde{\mathcal{H}}^{\otimes m_1} \otimes (\tilde{\mathcal{H}}^\vee)^{\otimes m_2} \otimes \bQ(m_3)_{\tilde{S}}\, ,
\]
which underlies a $\bQ$-VHS of weight $(m_1-m_2)k - 2m_3$ on~$\tilde{S}$. (The Tate twist $\bQ(m_3)_{\tilde{S}}$ is a constant local system; it only has an effect on how we index the Hodge filtration on $\tilde{\mathcal{T}}(\mathbf{m})_\bC$.) In particular, it makes sense to consider the set $Y(t) \subset \tilde{S}$ of points $x \in \tilde{S}$ for which the value~$t_x$ of~$t$ at~$x$ is a Hodge class. As remarked earlier, we can only have nonzero Hodge classes if the weight is zero, so we restrict our attention to triples $\mathbf{m} = (m_1,m_2,m_3)$ with $(m_1-m_2)k = 2m_3$. With this assumption on~$\mathbf{m}$, define
\[
Y(t) := \bigl\{x \in \tilde{S} \bigm| t_x \in \Fil^0 \tilde{\mathcal{T}}(\mathbf{m})_\bC(x) \bigr\}\, . 
\]
As $\Fil^0\tilde{\mathcal{T}}(\mathbf{m})_\bC \subset \tilde{\mathcal{T}}(\mathbf{m})_\bC$ is a holomorphic subbundle, $Y(t) \subset \tilde{S}$ is a countable union of closed irreducible analytic subsets of~$\tilde{S}$.

The first thing we deduce from this is that there is a countable union of proper analytic subsets $\Sigma \subset S$ such that the Mumford-Tate group $\MT_s$ is constant (in a suitable sense) on $S\setminus \Sigma$, and gets smaller if we specialize to a point in~$\Sigma$. Let us make this precise. Define $\tilde{\Sigma} \subset \tilde{S}$ by
\[
\tilde\Sigma := \cup\,  Y(t)
\]
where, with notation as introduced above, the union is taken over all triples $\mathbf{m}$ with $(m_1-m_2)k = 2m_3$, and all $t \in \Gamma(\tilde{S},\tilde{\mathcal{T}}(\mathbf{m}))$ such that $Y(t)$ is \emph{not} the whole~$\tilde{S}$. 

The main point of this definition of~$\tilde\Sigma$ is the following. If $x \in \tilde{S}\setminus \tilde{\Sigma}$ and $t_x$ is a Hodge class in the fiber of $\tilde{\mathcal{T}}(\mathbf{m})$ at the point~$x$, we can extend~$t_x$ in a unique way to a global section~$t$ of $\tilde{\mathcal{T}}(\mathbf{m})$ over~$\tilde{S}$, and by definition of~$\tilde{\Sigma}$ this global section is a Hodge class in every fiber. By contrast, if $x \in\tilde{\Sigma}$ then there is a tensor construction $\tilde{\mathcal{T}}(\mathbf{m})$ and a Hodge class $t_x$ whose horizontal extension~$t$ is not a Hodge class in every fiber. 

It follows that for $x \in \tilde{S}\setminus \tilde{\Sigma}$ the Mumford-Tate group $\MT_x \subset \GL(H) \times \Gm$ is independent of the choice of~$x$. We call the subgroup of $\GL(H) \times \Gm$ thus obtained the \emph{generic Mumford-Tate group}. Let us denote it by $\MT^\gen$. Further, for \emph{any} $x \in \tilde{S}$ we have an inclusion $\MT_x \subseteq \MT^\gen$.

The subset $\tilde{\Sigma} \subset \tilde{S}$ is stable under the action of the covering group of $\tilde{S}/S$ and therefore defines a subset $\Sigma \subset S$. It follows from the construction that $\Sigma$ is a countable union of closed analytic subspaces of~$S$. We say that a point $s \in S(\bC)$ is \emph{Hodge generic} if $s \notin \Sigma$. In somewhat informal terms we may restate the constancy of the Mumford-Tate group over $\tilde{S}\setminus \tilde{\Sigma}$ by saying that over Hodge generic locus the Mumford-Tate group is constant.

We can also draw conclusions pertaining to the loci in~$S$ where we have some given collection of Hodge classes. Start with a point $s_0 \in S(\bC)$, let $x_0 \in \tilde{S}$ be a point above~$s_0$, and consider a finite collection of nonzero classes~$t^{(i)}$, for $i=1,\ldots,r$, in tensor spaces $T(\mathbf{m}^{(i)}) = \Gamma\bigl(\tilde{S},\tilde{\mathcal{T}}(\mathbf{m}^{(i)})\bigr)$ that are Hodge classes for the Hodge structure at the point~$x_0$. With notation as above, the locus of points in~$\tilde{S}$ where all classes~$t^{(i)}$ are again Hodge classes is $Y(t^{(1)}) \cap \cdots \cap Y(t^{(r)})$. The image of this locus in~$S$ is a countable union of closed irreducible analytic subspaces. These components are called the Hodge loci of the given VHS.

\begin{definition}
\label{HodgeLociDef}
Let $\mathcal{H}$ be a polarizable $\bQ$-VHS over a complex manifold~$S$. A closed irreducible analytic subspace $Z \subseteq S$ is called a \emph{Hodge locus} of~$\mathcal{H}$ if there exist nonzero classes~$t^{(1)},\ldots,t^{(r)}$ in tensor spaces $T(\mathbf{m}^{(i)}) = \Gamma\bigl(\tilde{S},\tilde{\mathcal{T}}(\mathbf{m}^{(i)})\bigr)$ such that $Z$ is an irreducible component of the image of $Y(t^{(1)}) \cap \cdots \cap Y(t^{(r)})$ in~$S$.
\end{definition}

\begin{remark}
\label{CDKRemark}
If we start with a polarizable $\bZ$-VHS over a nonsingular complex algebraic variety~$S$ then by a theorem of Cattani, Deligne and Kaplan, see Corollary~1.3 in~\cite{Cat.Del.Kap}, the Hodge loci are algebraic subvarieties of~$S$. See also \cite{Voisin-HLoci}.
\end{remark}

We note that in the definition of Hodge loci there is no loss of generality to consider only a finite collection of Hodge classes. The requirement that a (possibly infinite) collection of classes~$t^{(i)}$ are all Hodge classes translates into the condition that the homomorphism $h \colon \bS \to \GL(H)_\bR$ that defines the Hodge structure factors through~$M_\bR$, where $M \subset \GL(H)$ is the common stabilizer of the given classes. Among the classes~$t^{(i)}$ we can then find a finite subcollection $t^{(i_1)},\ldots,t^{(i_n)}$ that have $M$ as their common stabilizer; so the Hodge locus we are considering is defined by these classes.

\section{Special subvarieties of~$\sfA_g$, and the Andr\'e-Oort conjecture}
\label{SpecSubvars}

In this section we discuss the notion of a special subvariety in a given Shimura variety. The abstract formalism of Shimura varieties provides a good framework for this but it has the disadvantage that it requires a lot of machinery. For this reason we shall give several equivalent definitions and we focus on concrete examples related to moduli of abelian varieties. In particular the version given in Definition~\ref{SpecSubvarDef3} can be understood without any prior knowledge of Shimura varieties.

Further we state the Andr\'e-Oort conjecture. A proof of this conjecture, under the assumption of the Generalized Riemann Hypothesis (GRH), has been announced by Klingler, Ullmo and Yafaev.  

\mosubsection{Shimura varieties}
\label{ShimVar}
In this paper, by a \emph{Shimura datum} we mean a pair $(G,X)$ consisting of an algebraic group~$G$ over~$\bQ$ together with a $G(\bR)$-conjugacy class~$X$ of homomorphisms $\bS \to G_\bR$, such that the conditions (2.1.1.1--3) of \cite{Del-VarSh} are satisfied. The space~$X$ is the disjoint union of finitely many connected components. These components have the structure of a hermitian symmetric domain of non-compact type.

Associated with such a datum is a subfield $E(G,X) \subset \bC$ of finite degree over~$\bQ$, called the reflex field. Given a compact open subgroup $K \subset G(\Afin)$ we have a Shimura variety $\Shim_K(G,X)$, which is a scheme of finite type over~$E(G,X)$, and for which we have
\[
\Shim_K(G,X)\bigl(\bC\bigr) = G(\bQ)\backslash\bigl(X \times G(\Afin)/K\bigr)\, .
\]
For $x \in X$ and $\gamma K \in G(\Afin)/K$ we denote by $[x,\gamma K] \in \Shim_K(G,X)\bigl(\bC\bigr)$ the class of $(x,\gamma K)$.

If $K^\prime \subset K$ is another compact open subgroup of~$G(\Afin)$ we have a natural morphism $\Shim_{K^\prime,K} \colon \Shim_{K^\prime}(G,X) \to \Shim_K(G,X)$. 

Let $\gamma \in G(\Afin)$. Given compact open subgroups $K$, $K^\prime \subset G(\Afin)$ with $K^\prime \subset \gamma K\gamma^{-1}$ we have a morphism $T_\gamma = [\cdot \gamma] \colon \Shim_{K^\prime}(G,X) \to \Shim_K(G,X)$ that is given on $\bC$-valued points by $T_\gamma [x,aK^\prime] = [x,a\gamma K]$. (For $\gamma=1$ we recover the morphism $\Shim_{K^\prime,K}$.) This induces a right action of the group $G(\Afin)$ on the projective limit
\[
\Shim(G,X) := \varprojlim_K \Shim_K(G,X)\, ,
\]
and $\Shim_K(G,X)$ can be recovered from $\Shim(G,X)$ as the quotient modulo~$K$. More generally, for compact open subgroups $K_1$, $K_2 \subset G(\Afin)$ and $\gamma \in G(\Afin)$ we have a Hecke correspondence~$T_\gamma$ from $\Shim_{K_1}(G,X)$ to $\Shim_{K_2}(G,X)$, given by the diagram
\begin{equation}
\label{TgDiag}
\Shim_{K_1}(G,X) \xleftarrow{\; \Shim_{K^\prime,K_1}\; } \Shim_{K^\prime}(G,X) \xrightarrow{\; [\cdot \gamma]\;} \Shim_{K_2}(G,X)\, ,
\end{equation}
where $K^\prime := K_1 \cap \gamma K_2\gamma^{-1}$.

\mosubsection{The description of~$\cA_g$ as a Shimura variety}
\label{AgShimVar}
With notation as in~\ref{Notation}, let $\Siegel_g$ denote the space of homomorphisms $h \colon \bS \to \CSp_{2g,\bR}$ that define a Hodge structure of type $(-1,0) + (0,-1)$ on~$V$ for which $\pm (2\pi i)\cdot \phi \colon V \times V \to \bZ(1)$ is a polarization. The group $\CSp_{2g}(\bR)$ acts transitively on~$\Siegel_g$ by conjugation, and the pair $(\CSp_{2g,\bQ},\Siegel_g)$ is an example of a Shimura datum. The reflex field of this datum is~$\bQ$.

The associated Shimura variety is known as the \emph{Siegel modular variety} and may be identified with the moduli space of principally polarized abelian varieties with a level structure. To explain this in detail, we define, for a positive integer~$m$, a compact open subgroup $K_m \subset G(\Afin)$ by 
\[
K_m = \bigl\{\gamma \in \CSp(V \otimes \Zhat,\phi) \bigm| \gamma \equiv \id \pmod{m}\bigr\}\, .
\]
For $m \geq 3$ we then have an isomorphism 
\[
\beta\colon \mathcal{A}_{g,[m],\bQ} \isomarrow \Shim_{K_m}(\CSp_{2g,\bQ},\Siegel_g)\, .
\]
On $\bC$-valued points $\beta$ is given as follows. To begin with, a $\bC$-valued point of~$\mathcal{A}_{g,[m]}$ is a triple $(A,\lambda,\alpha)$ consisting of a complex ppav of dimension~$g$ together with a level~$m$ structure that is symplectic for some choice of a primitive $m$th root of unity $\zeta \in \bC$. Write $H := H_1(A,\bZ)$ for the singular homology group in degree~$1$ of (the complex manifold associated with)~$A$, and let $\psi \colon H \times H \to \bZ(1)$ be the polarization associated with~$\lambda$; see~(\ref{AVEqCat}). Choose a symplectic similitude $s \colon H \isomarrow V$. (Even though $\psi$ takes values in~$\bZ(1)$ and $\phi$ takes values in~$\bZ$, the notion of a symplectic similitude makes sense.) Via~$s$, the natural Hodge structure on~$H$ corresponds to an element $x \in \Siegel_g$. Further, we have an identification $A[m] \isomarrow H/mH$, such that $e^\lambda(P,Q) = \exp({\psi(y,z)/m})$ if $P$, $Q \in A[m]\bigl(\bC\bigr)$ correspond to $y$, $z \in H/mH$, respectively. (Note that $\exp({\psi(y,z)/m})$ is well-defined.) Hence the given level structure~$\alpha$ can be viewed, via~$s$, as an element $\gamma K_m \in \CSp(V \otimes (\bZ/m\bZ),\phi) = \CSp(V \otimes \Zhat,\phi)/K_m$. The isomorphism~$\beta$ then sends $(A,\lambda)$ to the point $[x,\gamma K_m]$.

For $m \geq 3$ the group $\CSp_{2g}(\bQ)$ acts properly discontinuously on $\Siegel_g \times \CSp_{2g}(\Afin)/K_m$. The space $\Siegel_g$ has two connected components, and if $\Siegel_g^+$ is one of these, the $\varphi(m)$ components of $\mathcal{A}_{g,[m],\bC}$ all have the form $\Gamma\backslash \Siegel_g^+$ for an arithmetic subgroup $\Gamma \subset \CSp_{2g}(\bQ)$. For $m=1$ it is no longer true that $\CSp_{2g}(\bQ)$ acts properly discontinuously on $\Siegel_g \times \CSp_{2g}(\Afin)/K_1$, and we should interpret the quotient space $\CSp_{2g}(\bQ)\backslash\bigl(\Siegel_g \times \CSp_{2g}(\Afin)/K_1\bigr)$ as an orbifold. Alternatively, we may take the actual quotient space; this gives us an isomorphism
\[
\CSp_{2g}(\bQ)\backslash\bigl(\Siegel_g \times \CSp_{2g}(\Afin)/K_1\bigr) \isomarrow \sfA_g(\bC)
\]
between $\Shim_{K_1}(\CSp_{2g,\bQ},\Siegel_g)\bigl(\bC\bigr)$ and the set of $\bC$-valued points of the coarse moduli space.

\begin{remark}
\label{BorelEmb}
Let $\mathcal{L} = \mathcal{L}(V_\bQ,\phi)$ denote the symplectic Grassmanian of Lagrangian (i.e., maximal isotropic) subspaces of~$V_\bQ$ with respect to the form~$\phi$. The map $\Siegel_g \to \mathcal{L}(\bC)$ that sends a Hodge structure $y \in \Siegel_g$ to the corresponding Hodge filtration $\Fil^0 \subset V_\bC$ is an open immersion, known as the Borel embedding. For an arbitrary Shimura datum $(G,X)$ we have, in a similar manner, a Borel embedding of~$X$ into the $\bC$-points of a homogeneous projective variety.

Because $\Fil^0 \subset V_\bC$ is maximal isotropic, the form~$\phi$ induces a perfect pairing $\bar{\phi} \colon \Fil^0 \times V_\bC/\Fil^0 \to \bC$. Via the Borel embedding, the tangent space of $\Siegel_g$ at the point~$y$ maps isomorphically to the tangent space of~$\mathcal{L}$ at the point~$\Fil^0$, which gives us
\begin{equation}
\label{TySiegel}
\begin{split}
& T_y(\Siegel_g) \isomarrow {}\Hom^\sym(\Fil^0,V_\bC/\Fil^0)\\ &{}\qquad := \bigl\{\beta \colon \Fil^0 \to V_\bC/\Fil^0 \bigm| \bar{\phi}(v,\beta(v^\prime)) = \bar{\phi}(v^\prime,\beta(v))\  \mbox{for all $v$, $v^\prime \in \Fil^0$} \bigr\}\, .
\end{split}
\end{equation}
The condition that $\bar{\phi}(v,\beta(v^\prime)) = \bar{\phi}(v^\prime,\beta(v))$ means that $\beta$ is its own dual, via the isomorphisms $\Fil^0 \isomarrow (V_\bC/\Fil^0)^\vee$ and $V_\bC/\Fil^0 \isomarrow (\Fil^0)^\vee$ induced by~$\bar{\phi}$; whence the notation $\Hom^\sym$.
\end{remark}

\mosubsection{Special subvarieties}
\label{SpecialSubvars}
There are several possible approaches to the notion of a special subvariety of a given Shimura variety. Let us start with the construction that gives the quickest definition. We use the language of Shimura varieties; the example to keep in mind is the Siegel modular variety $\Shim(\CSp_{2g,\bQ},\Siegel_g)$. After this we shall give some alternative definitions that lean less heavily on the abstract formalism of Shimura varieties. The reader who prefers to avoid the language of Shimura varieties is encouraged to skip ahead to Definition~\ref{SpecSubvarDef3}.

We fix a Shimura datum $(G,X)$ and a compact open subgroup $K \subset G(\Afin)$. If $(M,Y)$ is a second Shimura datum, by a morphism $f \colon (M,Y) \to (G,X)$ we mean a homomorphism of algebraic groups $f \colon M \to G$ such that for any $y \colon \bS \to M_\bR$ in~$Y$ the composite $f \circ y \colon \bS \to G_\bR$ is an element of~$X$. The existence of such a morphism of Shimura data implies that $E(M,Y)$ contains $E(G,X)$, and $f$ gives rise to a morphism of schemes
\[
\Shim(f) \colon \Shim(M,Y) \to \Shim(G,X)_{E(M,Y)}\, .
\]

Let $\gamma \in G(\Afin)$. The image of the composite morphism
\begin{equation}\label{TgSh(f)}
\Shim(M,Y)_\bC \xrightarrow{\;\Shim(f)\;} \Shim(G,X)_\bC \xrightarrow{\; [\cdot \gamma]\;} \Shim(G,X)_\bC \longrightarrow \Shim_K(G,X)_\bC
\end{equation}
is a reduced closed subscheme of $\Shim_K(G,X)_\bC$.

\begin{definition}[Version~1]
\label{SpecSubvarDef1}
A closed subvariety $Z \subset \Shim_K(G,X)_\bC$ is called a \emph{special subvariety} if there exists a morphism of Shimura varieties $f \colon (M,Y) \to (G,X)$ and an element $\gamma \in G(\Afin)$ such that $Z$ is an irreducible component of the image of the morphism~\eqref{TgSh(f)}. 
\end{definition}

For a second approach we fix an algebraic subgroup $M \subset G$ over~$\bQ$. Define $Y_M \subset X$ as the set of all $x \colon \bS \to G_\bR$ in~$X$ that factor through $M_\bR \subset G_\bR$. We give $Y_M$ the topology induced by the natural topology on~$X$. The group $M(\bR)$ acts on~$Y_M$ by conjugation. It can be shown (see \cite{Moonen-PhD}, Section~I.3, or \cite{Moonen-LinI}, 2.4) that $Y_M$ is a finite union of orbits under~$M(\bR)$. We remark that the condition that $Y_M$ is nonempty imposes strong restrictions on~$M$; it implies, for instance, that $M$ is reductive.

\begin{definition}[Version~2]
\label{SpecSubvarDef2}
A closed subvariety $Z \subset \Shim_K(G,X)_\bC$ is called a \emph{special subvariety} if there exists an algebraic subgroup $M \subset G$ over~$\bQ$, a connected component $Y^+ \subset Y_M$, and an element $\gamma \in G(\Afin)$ such that $Z(\bC) \subset \Shim_K(G,X)\bigl(\bC\bigr)$ is the image of $Y^+ \times \{\gamma K\} \subset X \times G(\Afin)/K$ under the natural map to $\Shim_K(G,X)\bigl(\bC\bigr) = G(\bQ)\backslash \bigl(X \times G(\Afin)/K \bigr)$.
\end{definition}

For the equivalence of this definition with Version~1 we refer to \cite{Moonen-PhD}, Proposition~3.12 or \cite{Moonen-LinI}, Remark~2.6.

Our third version of the definition, which in some sense is the most conceptual one, describes special subvarieties of $\Shim_K(G,X)$ as the Hodge loci of certain natural VHS associated with representations of the group~$G$. As we wish to highlight the Siegel modular case, we only state this version of the definition for special subvarieties of~$\mathcal{A}_{g,[m],\bC}$ for some $m \geq 3$. (See, however, Remark~\ref{Version3Rem}.)

In order to talk about Hodge loci we need some VHS. The most natural $\bQ$-VHS on $\mathcal{A}_{g,[m],\bC}$ to consider is the variation~$\mathcal{H}$ whose fiber at a point $(A,\lambda,\alpha)$ is given by $H_1(A,\bQ)$. (One could also work with the $\bQ$-VHS given by the first cohomology groups; for what we want to explain it makes no difference.) The Hodge loci of this VHS are precisely the special subvarieties.

\begin{definition}[Version~3]
\label{SpecSubvarDef3}
A closed subvariety $Z \subset \mathcal{A}_{g,[m],\bC}$ is called a \emph{special subvariety} if it is a Hodge locus of the $\bQ$-VHS $\mathcal{H}$ over $\mathcal{A}_{g,[m],\bC}$ whose fiber at a point $(A,\lambda,\alpha)$ is $H_1(A,\bQ)$.
\end{definition}

Concretely this means that the special subvarieties are ``defined by'' the existence of certain Hodge classes, i.e., they are the maximal closed irreducible subvarieties of $\mathcal{A}_{g,[m],\bC}$ on which certain given classes are Hodge classes. In order to make this precise, one has to pass to the universal cover, where the underlying local system can be trivialized, as explained in Section~\ref{HodgeLoci}. Note that in this case the Hodge loci are algebraic subvarieties of~$\mathcal{A}_{g,[m],\bC}$; this can be shown by proving that this notion of a special subvariety agrees with the one given by Version~1 of the definition but it also follows from the theorem of Cattani, Deligne and Kaplan mentioned in Remark~\ref{CDKRemark}.

\begin{remark}
\label{Version3Rem}
Version~3 of the definition of a special subvariety, in terms of Hodge loci, can be extended without much difficulty to arbitrary Shimura varieties $\Shim_K(G,X)$. The variations of Hodge structure one considers are those associated with representations of the group~$G$. See for instance \cite{Moonen-LinI}, Proposition~2.8, which also proves the equivalence with the earlier definitions. We note that, even if one is only interested in statements about special subvarieties of~$\mathcal{A}_g$, there are good reasons to extend this notion to more general Shimura varieties. The formalism of Shimura varieties enables one to perform some useful constructions, such as the passage to an adjoint Shimura datum, or the reduction of a problem to the case of a simple Shimura datum.
\end{remark}

\begin{example}
\label{SpecSubvarEndom}
Let $(A,\lambda,\alpha)$ be a complex ppav of dimension~$g$ with a symplectic level~$m$ structure, for some $m \geq 3$. Let $D := \End^0(A)$ be its endomorphism algebra. We should like to describe the largest closed (irreducible) subvariety $Z \subset \mathcal{A}_{g,[m],\bC}$ that contains the moduli point $[A,\lambda,\alpha] \in \mathcal{A}_{g,[m]}\bigl(\bC\bigr)$, and such that all endomorphisms of~$A$ extend to endomorphisms of the universal abelian scheme over~$Z$. In terms of Hodge classes this means that all elements of~$D$, viewed as Hodge classes in $\End_\bQ\bigl(H_1(A,\bQ)\bigr)$, should extend to Hodge classes in the whole $\bQ$-VHS over~$Z$ given by the endomorphisms of the first homology groups. (Note that we may pass to $\bQ$-coefficients, for if $f \in \End(A)$ is an endomorphism such that some positive multiple $nf$ extends to an endomorphism of the whole family, $f$ itself extends.)

As before, choose a symplectic similitude $s\colon H_1(A,\bZ) \isomarrow V$ with respect to the polarization forms on both sides. This gives~$V_\bQ$ the structure of a left module over the algebra~$D$. Let $h_0 \in \Siegel_g$ be the point given by the Hodge structure on $H_1(A,\bZ)$, viewed as a Hodge structure on~$V$ via~$s$. Further, let $M \subset \CSp_{2g,\bQ} = \CSp(V_\bQ,\phi)$ be the algebraic subgroup given by
\[
M = \CSp(V_\bQ,\phi) \cap \GL_D(V_\bQ)\, .
\]
The homomorphism $h_0 \colon \bS \to \CSp_{2g,\bR}$ factors through~$M_\bR$. Conversely, if $h \in \Siegel_g$ factors through~$M_\bR$ then $D$ acts by endomorphisms on the corresponding abelian variety. This means that $Z \subset \Shim_{K_m}(\CSp_{2g,\bQ},\Siegel_g)_\bC = \mathcal{A}_{g,[m],\bC}$ is the special subvariety that is obtained, with notation as in Definition~\ref{SpecSubvarDef2}, as the image of $Y_M^+ \times \{\gamma K_m\}$ in $\mathcal{A}_{g,[m],\bC}$, were $Y_M^+ \subset Y_M$ is the connected component containing~$h_0$ and $\gamma K_m \in \CSp_{2g}(\Afin)/K_m$ is the class corresponding (via~$s$) with the given level structure~$\alpha$.

We conclude that the closed subvarieties $Z \subset \mathcal{A}_{g,[m],\bC}$ ``defined by'' the existence of endomorphisms, as made precise above, are examples of special subvarieties.
\end{example}

\begin{remark}
\label{PELData}
The special subvarieties considered in Example~\ref{SpecSubvarEndom} are referred to as special subvarieties \emph{of PEL type}; the name comes from the fact that they have a modular interpretation in terms of abelian varieties with a \emph{p}olarization, given \emph{e}ndomorphisms and a \emph{l}evel structure.

In the above example our focus is on the concrete modular interpretation of these special subvarieties. For many purposes it is relevant to also have a good description of these examples in the language of Shimura varieties; see also \cite{Del-TravShim}, Section~4. Here one starts with data $(D,*,V,\phi)$ where $D$ is a finite dimensional semisimple $\bQ$-algebra, $*$ is a positive involution, $V$ is a faithful (left) $D$-module of finite type, and $\phi \colon V \times V \to \bQ$ is an alternating form such that 
\begin{equation}
\label{phidvv}
\phi(dv,v^\prime) = \phi(v,d^* v^\prime)\qquad \mbox{for all $v$, $v^\prime \in V$ and $d \in D$.}
\end{equation}
Let $G = \CSp(V,\phi) \cap \GL_D(V)$ be the group of $D$-linear symplectic similitudes of~$V$. Next one considers a $G^0(\bR)$-conjugacy class~$X$ of homomorphisms $\bS \to G_\bR$ defining on~$V$ a Hodge structure of type $(-1,0) + (0,-1)$ such that $\pm 2\pi i\cdot \phi$ is a polarization. The pair $(G^0,X)$ is then a Shimura datum.

In Example~\ref{SpecSubvarEndom}, the data $(D,V,\phi)$ are given, and we take for~$*$ the involution on~$D$ defined by~(\ref{phidvv}). For~$X$ we take the conjugacy class of the homomorphism~$h_0$. By construction, the inclusion $G \hookrightarrow \CSp_{2g}$ defines a morphism of Shimura data $f\colon (G^0,X) \to (\CSp_{2g},\Siegel_g)$, and the special subvariety $Z \subset \mathcal{A}_{g,[m],\bC}$ of~\ref{SpecSubvarEndom} is an irreducible component of a Hecke translate of the image of $\Shim(G^0,X)$ in $\Shim_{K_m}(\CSp_{2g,\bQ},\Siegel_g)_\bC = \mathcal{A}_{g,[m],\bC}$.
\end{remark}

\begin{example}
Let $(G,X)$ be a Shimura datum. If $x \colon \bS \to G_\bR$ is an element of~$X$, we define the Mumford-Tate group of~$x$ to be the smallest algebraic subgroup $\MT_x \subset G$ over~$\bQ$ such that $x$ factors through~$\MT_{x,\bR}$.

Let $s = [x,\gamma K] \in \Shim_K(G,X)\bigl(\bC\bigr)$. Up to conjugation by an element of~$G(\bQ)$, the Mumford-Tate group $\MT_x \subset G$ is independent of the choice of the chosen representative $(x,\gamma K)$ for~$s$. Hence we may define the Mumford-Tate group of~$s$ to be $\MT_s := \MT_x$. 

A point $x \in X$ is called a \emph{special point} if $\MT_x$ is a torus. A point $s \in \Shim_K(G,X)\bigl(\bC\bigr)$ is a special subvariety of $\Shim_K(G,X)_\bC$ if and only if $\MT_s$ is a torus, i.e., if for some (equivalently, any) representative $[x,\gamma K]$ the point $x \in X$ is special. 

In the Siegel modular variety $\sfA_g$ over~$\bC$, the special points are precisely the CM points, i.e., the points corresponding to ppav $(A,\lambda)$ such that $A$ is an abelian variety of CM type; see Mumford \cite{Mumf-ShimPaper}, \S~2.
\end{example}

\mosubsection{Basic properties of special subvarieties}
\label{BasicPropSpecial}
We list a number of elementary properties.
\smallskip

\noindent
(a) \emph{Hecke images of special subvarieties are again special.} More precisely, consider a Hecke correspondence~$T_\gamma$ as given by~\eqref{TgDiag}. Let $Z \subset \Shim_{K_1}(G,X)_\bC$ be a special subvariety. Write $T_\gamma(Z)$ for the image of $\Shim_{K^\prime,K_1}^{-1}(Z)$ under the map $[\cdot \gamma]$. Then all irreducible components of $T_\gamma(Z)$ are special subvarieties of $\Shim_{K_2}(G,X)_\bC$.

As a particular case, if we have compact open subgroups $K^\prime \subset K \subset G(\Afin)$ and if $Y \subset \Shim_{K^\prime}(G,X)_\bC$ is a special subvariety, the image of~$Y$ in $\Shim_K(G,X)_\bC$ is again a special subvariety. Conversely, if $Z \subset \Shim_K(G,X)_\bC$ is a special subvariety, the irreducible components of $\Shim_{K^\prime,K}^{-1}(Z)$ are special subvarieties of $\Shim_{K^\prime}(G,X)_\bC$. In the study of special subvarieties, these remarks often allow us to choose the level subgroup~$K$ as small as needed.
\smallskip

\noindent
(b) \emph{The special points in a special subvariety are dense.} If $Z \subset \Shim_K(G,X)_\bC$ is a special subvariety then the special points in~$Z$ are dense for the analytic topology on~$Z(\bC)$; in particular they are Zariski dense. In order to prove this, the essential point is to show that $Z$ contains at least \emph{one} special point. The density of special points then follows from the fact that, with notation as in~Def.~\ref{SpecSubvarDef2}, the set of special points in~$Y_M$ is stable under the action of $M(\bQ)$, and that $M(\bQ)$ is analytically dense in~$M(\bR)$.

In view of the importance of special points, let us sketch a proof of the existence of at least one special point, following \cite{Mumf-ShimPaper}, \S~3. The argument uses the fact from the theory of reductive groups that, given a maximal torus $T \subset G_\bR$, the $G(\bR)$-conjugacy class of~$T$ contains a maximal torus that is defined over~$\bQ$. 

Start with any $x \colon \bS \to G_\bR$ in~$X$. Let $C \subset G_\bR$ be the centralizer of $x(\bS)$, which is a connected reductive subgroup of~$G_\bR$. (See \cite{Springer-LAG}, Lemma 15.3.2.) Choose a maximal torus $T \subset C$. Because $x(\bS)$ is contained in the center of~$C$, we have $x(\bS) \subseteq T$. Hence, if $T^\prime \subset G_\bR$ is any torus that contains~$T$ then $T^\prime$ centralizes~$x(\bS)$ and therefore $T^\prime \subset C$. It follows that $T$ is also a maximal torus of~$G_\bR$. By the general fact stated above, there exists an element $g \in G(\bR)$ and a maximal torus $S \subset G$ (over~$\bQ$) such that $gTg^{-1} = S_\bR$. Then $gx = \Inn(g) \circ x \colon \bS \to G_\bR$ factors through $S_\bR$; hence $gx \in X$ is a special point.

A more refined version of the existence of special points plays a role in the theory of canonical models of Shimura varieties. See \cite{Del-TravShim}, Theorem~5.1.
\smallskip

\noindent
(c) \emph{Intersections of special subvarieties are again special.} This is easily seen using Version~2 of the definition.
\smallskip

\noindent
(d) \emph{After passage to an appropriate level cover, special subvarieties are locally symmetric.} Suppose $Z \subset \Shim_K(G,X)$ is a special subvariety. We may then find a subgroup $K^\prime \subset K$ of finite index and an irreducible component $Z^\prime \subset \Shim_{K^\prime}(G,X)$ of the inverse image of~$Z$ such that $Z^\prime$ is a locally symmetric variety. (We shall recall the definition of this notion in \ref{HainResults} below.)

\begin{remark}
\label{phiTSpecRem}
A situation we often encounter is that we have an abelian scheme $A \to T$ over some complex algebraic variety~$T$ (assumed to be irreducible), with a principal polarization $\lambda \colon A\to A^t$ and a level~$m$ structure~$\alpha$. This gives us a morphism $\tau \colon T \to \mathcal{A}_{g,[m],\bC}$, and we should like to know if the scheme-theoretic image of~$\tau$ is a special subvariety. (In this situation, the scheme-theoretic image is the reduced closed subscheme of $\mathcal{A}_{g,[m]}$ that has as underlying set the Zariski closure of the topological image of~$\tau$.) 

We remark that the answer to this question only depends on~$A$ up to isogeny, and is independent of the polarization and the level structure. 

Of course, if we change the polarization or the level structure, or if we replace $A/T$ by an isogenous abelian scheme, the morphism~$\tau$ is replaced by another morphism $\tau^\prime\colon T \to \mathcal{A}_{g,[m]}$, but if the image of~$\tau$ is special, so is the image of~$\tau^\prime$. (For simplicity of exposition, we here assume the new polarization is again principal, but even this assumption can be dropped.) The reason that this is true is that the Mumford-Tate group of an abelian variety~$A$ only depends on $A$ up to isogeny, and not on any additional structures. 

The conclusion, then, is that one may set up the situation as one finds it convenient. We could equip $(A,\lambda)$ with a level structure (possibly after replacing~$T$ with a cover), allowing us to work with a fine moduli scheme $\mathcal{A}_{g,[m]}$. Alternatively, one could forget about the level structure and consider the morphism $T \to \mathcal{A}_g$ to the moduli stack; in this case we should talk about special substacks of~$\mathcal{A}_g$, which makes perfectly good sense. Yet another option is to consider the morphism $T \to \sfA_g$ to the coarse moduli scheme. For the question whether the closed image of~$T$ is special, it does not matter which version we consider.

Further, since we are only interested in the closure of $\tau(T)$, we may replace~$T$ by an open subset, and in fact it suffices to know the generic fiber of~$A/T$.
\end{remark}

\begin{conjecture}[Andr\'e-Oort, \cite{Andre-1989}, \cite{FO-Que}, \cite{FO-Cortona}]
\label{AOConj}
Let $(G,X)$ be a Shimura datum. Let $K \subset G(\Afin)$ be a compact open subgroup and let $Z \subset \Shim_K(G,X)_\bC$ be an irreducible closed algebraic subvariety such that the special points on~$Z$ are dense for the Zariski topology. Then $Z$ is a special subvariety.
\end{conjecture}

An alternative way of stating the conjecture is that, given a set of special points $\mathfrak{S} \subset \Shim_K(G,X)\bigl(\bC\bigr)$, the irreducible components of the Zariski closure of~$\mathfrak{S}$ are special subvarieties.

We refer to Noot's Bourbaki lecture~\cite{Noot-Bourbaki} for an excellent overview of what was known about the conjecture in 2004.

Klingler and Yafaev~\cite{Klingler.Yaf}, using the work of Ullmo and Yafaev~\cite{Ullmo.Yaf}, have announced a proof of the Andr\'e-Oort  conjecture, assuming the Generalized Riemann hypothesis (GRH) for CM fields.

\begin{theorem}[Klingler-Yafaev]
\label{KlingYafThm}
Assume the GRH holds for all CM fields. Then the Andr\'e-Oort conjecture is true.
\end{theorem}

There are some special cases of the Andr\'e-Oort conjecture that are known to hold without any assumptions on the GRH. See for instance \cite{Moonen-LinII}, Theorem~5.7, later refined by Yafaev in~\cite{Yafaev-Moonen}, Theorem~1.2 of \cite{Edix.Yaf}, and Theorem~1.2.1 (with condition~(2)) of \cite{Klingler.Yaf}. In all these cases, however, further assumptions are needed on the set of special points. Apart from trivial cases, the only completely unconditional case of the Andr\'e-Oort conjecture that we know of, is the main result of Andr\'e's paper~\cite{Andre-pairs}, which proves the conjecture for subvarieties of a product of two modular curves, and the more recent extension of this by Pila~\cite{Pila} to arbitrary products of modular curves.

\begin{remark}
\label{A1xA2Rem}
As in Remark~\ref{phiTSpecRem}, consider a principally polarized abelian scheme $(A,\lambda)$ over some complex algebraic variety~$T$. Again we consider the question whether the closed image of the morphism $\tau \colon T \to \sfA_g$ is a special subvariety. Suppose that $A/T$ is isogenous to a product $A_1 \times A_2$ with $\dim(A_i/T) = g_i$. After choosing polarizations (principal, say) on the factors~$A_i$ we get morphisms $\tau_i \colon T \to \sfA_{g_i}$. In this situation, if the closure of $\tau(T)$ in~$\sfA_g$ is a special subvariety, the closure of $\tau_i(T)$ in $\sfA_{g_i}$ is special, too. If one believes the Andr\'e-Oort conjecture this is clear, and it is in fact not very difficult to show this using our definitions of a special subvariety. 

The converse implication does not hold, in general. So, if the closures of $\tau_1(T)$ and $\tau_2(T)$ are special, this does not imply that the closure of $\tau(T)$ is special. Looking at it from the perspective of the Andr\'e-Oort conjecture, suppose that for each of the two factors $A_i/T$ separately, we have a Zariski dense collection of points in~$T(\bC)$ at which the fiber of $A_i/T$ is of CM type. Then it is not true, in general, that there is a dense set of points at which the fibers are \emph{simultaneously} of CM type. (For a concrete example, take $T$ to be an open part of $\bA^1\setminus\{0,1\}$, let $A_1/T$ be a family of elliptic curves that is not isotrivial, and take for $A_2$ the same family with a suitable shift in the parameter, i.e., such that $A_{2,t} = A_{1,t+\epsilon}$ for some fixed $\epsilon \in \bC$.)
\end{remark}

\section{Special subvarieties in the Torelli locus}
\label{SsT}

Throughout this section we work over~$\bC$.

\begin{conjecture}[Coleman, \cite{Col}, Conjecture~6]
\label{ColemanConj}
Given $g \geq 4$ there are only finitely many non-singular projective curves~$C$ over~$\bC$, up to isomorphism, of genus  $g$ and the Jacobian~$J_C$ is a CM abelian variety.
\end{conjecture}

As we shall discuss below, the conjecture is known to be false for $g \leq 7$, so a corrected version of the conjecture should have as an assumption that $g \geq 8$.

Both Coleman's conjecture~\ref{ColemanConj} and the Andr\'e-Oort conjecture~\ref{AOConj} were inspired by the analogy with the Manin-Mumford conjecture, now a theorem of Raynaud. For proofs of the Manin-Mumford conjecture, see \cite{Raynaud1}, \cite{Raynaud2}; for an easy proof see \cite{Pink.R}. The analogy between this conjecture and the Andr\'e-Oort conjecture is discussed for instance in \cite{Moonen-Durham}, Section~6 and in \cite{Noot-Bourbaki}, Section~3. As we shall discuss next, the Andr\'e-Oort conjecture also has important implications for Coleman's conjecture.

\begin{expectation}[Oort, \cite{FO-Cortona}, \S~5]
\label{ZinTorelli}
For large~$g$ (in any case $g \geq 8$), there does not exist a special subvariety $Z \subset \sfA_{g}$ with $\dim(Z) \geq 1$ such that $Z \subseteq \Torelli_g$ and $Z \cap \Torelli_g^\open$ is nonempty.
\end{expectation}

Note that the assumption that $Z \cap \Torelli_g^\open$ is nonempty implies that this intersection is open and dense in~$Z$.

\begin{remark}
The condition that $Z$ meets $\Torelli_g^\open$ is important. Indeed, it is easy to see that for any $g \geq 2$ there exist special subvarieties $Z \subset \sfA_g$ of positive dimension with $Z \subset \Torelli_g$. For $g \leq 3$ this is clear, as $\Torelli_g = \sfA_g$ is special. Assume $g > 3$. Choose a base variety~$T$ and a stable curve $C \to T$ of genus~$3$, such that the closure of the image of~$T$ in~$\sfA_3$ is a special subvariety of positive dimension~$d$. Let $J\to T$ be the Jacobian, $\lambda$ its canonical principal polarization. Next take an elliptic curve $E$ with complex multiplication, and let $\mu$ be its principal polarization. Then for every $t \in T(\bC)$ the moduli point $\xi_t \in \sfA_g(\bC)$ of the ppav $(J_t,\lambda) \times (E,\mu)^{g-3}$ lies in the closed Torelli locus~$\Torelli_g$, as it is the Jacobian of the curve that is obtained from~$C_t$ by attaching to it a tail of $g-3$ copies of~$E$. Moreover, it is not hard to see that the closure of the set of points~$\xi_t$ is a $d$-dimensional special subvariety of~$\sfA_g$. In this way we can produce many positive dimensional special subvarieties in $\Torelli_g$ for any $g\geq 2$. 
\end{remark}

\begin{remark}  
Assume we have an integer~$g$ for which Expectation~\ref{ZinTorelli} holds. Assume furthermore that the Andr\'e-Oort Conjecture~\ref{AOConj} is true.  Then Coleman's Conjecture~\ref{ColemanConj} is true in genus~$g$. In fact, let $\mathrm{CM}(\Torelli_g^\open) \subset \Torelli(\bC)$ be the set of all CM Jacobians of dimension~$g$. If this set is infinite, its Zariski closure in~$\sfA_g$ has at least one irreducible component $Z$ of positive dimension, which by~\ref{AOConj} is special, contradicting~\ref{ZinTorelli}. Hence \ref{ZinTorelli} and~\ref{AOConj} together imply~\ref{ColemanConj}.
\end{remark}

\begin{remark}
For $g>3$ we know that $\Torelli_g$ itself is not a special subvariety, and in fact there is no special subvariety $S \subsetneq \sfA_g$ that contains~$\Torelli_g$. (Note that $\Torelli_g \neq \sfA_g$ for $g>3$.) The simplest argument we know for this is to use information about the geometric monodromy. For convenience, let us pass to moduli spaces with a level~$m$ structure, for some $m \geq 3$. Write $\cT_{g,[m]} \subset \cA_{g,[m]}$ for the Torelli locus in~$\cA_{g,[m]}$. Over~$\cA_{g,[m]}$ we have the natural $\bQ$-VHS~$\cH$ considered before (see just before Definition~\ref{SpecSubvarDef3}), and the assertion that $\Torelli_g$ is not contained in any special subvariety $S \subsetneq \sfA_g$ is implied by the fact that the generic Mumford-Tate group of the restriction of~$\cH$ to~$\cT_{g,[m]}$ is the full group $\CSp_{2g,\bQ}$. 

Write $\cH^\prime$ for the restriction of~$\cH$ to~$\cT_{g,[m]}$. Choose a Hodge generic base point $b \in \cT_{g,[m]}$, let $H$ be the fiber of~$\cH^\prime$ at~$b$, and let $\psi$ be the polarization form on~$H$. Further, let $M \subset \CSp_{2g,\bQ} = \CSp(H,\psi)$ be the Mumford-Tate group at~$b$, which is the generic Mumford-Tate group of~$\cH^\prime$. Then $M$ contains $\bG_\mult \cdot 1$. On the other hand we have a monodromy representation
\[
\rho \colon \pi_1(\cT_{g,[m]},b) \to \GL(H)\, ,
\]
and since we have passed to a level cover, $\rho$ factors through $\Sp(H,\psi)$. By a result of Deligne, see \cite{Del-WeilK3}, Proposition~7.5, the image of~$\rho$ has a subgroup of finite index that is contained in~$M(\bC)$. So we are done if we can show that the image of~$\rho$ is Zariski dense in $\Sp(H,\psi)$. This can be seen, for instance, using transcendental methods. In fact, the homomorphism $\rho \colon \pi_1(\cT_{g,[m]},b) \to \Sp(H,\psi)$ may be identified with the natural homomorphism $\Gamma_g \to \Sp(H,\psi)$ from the mapping class group in genus~$g$ to the symplectic group, and it is a classical result that this homomorphism is surjective. (See for instance \cite{ACG-GAC2}, Chap.~15, especially~\S~3.)

The argument given here is entirely based on information about the geometric monodromy, and there are other loci in~$\Torelli_g$ to which the same reasoning applies. As an example, if $\HElocus_g \subset \Torelli_g$ is the hyperelliptic locus, the image of the geometric monodromy representation on~$\HElocus_g$ is again dense in the symplectic group $\Sp(H,\psi)$; see \cite{ACampo}, Theorem~1. By the above argument, it follows that there is no special subvariety $S \subsetneq \sfA_g$ that contains~$\HElocus_g$.
\end{remark}

\mosubsection{A modified version of Coleman's Conjecture}
\label{ModifColeman}
In \cite{Chai.FO-AVIsogJ} we find a modified version of \ref{ColemanConj} which does hold, at least conditionally, for every $g>3$. A $g$-dimensional abelian variety $A$, say over~$\bC$, is said to be a \emph{Weyl CM abelian variety} if $L := \End^0(A)$ is a field of degree $2g$ over~$\bQ$ whose Galois closure has degree $2^g\cdot g!$ over~$\bQ$. It can be shown that, in a suitable sense, most CM abelian varieties are of this type.

\begin{theorem}[Chai and Oort]
\label{WeylCM}
Assume the Andr\'e-Oort Conjecture~{\upshape{\ref{AOConj}}} to be true. Then for every $g>3$ the number of  Weyl CM points in $\Torelli_g^\open$ is finite.
\end{theorem}

See \cite{Chai.FO-AVIsogJ}, 3.7. We remark that for $g \geq 4$ we do not know any example of a Jacobian of Weyl CM type. In connection with this, note that for a non-hyperelliptic curve~$C$ of genus $g(C) > 1$ with $\Aut(C)\neq \{\id\}$ the Jacobian $J_C$ does not give a Weyl CM point.

\mosubsection{Results of Hain}
\label{HainResults}
In his paper \cite{Hain-Locsymm}, Hain proved some results inspired by~\ref{ZinTorelli}. Though the results are conditional, they point in an interesting direction, as they suggest that ball quotients should play a special role. (The open unit ball in~$\bC^n$ is the symmetric space $\SU(n,1)/\UU(n)$.) It should be pointed out that all non-trivial examples known to us, in genus at least~$4$, are indeed ball quotients; see the next sections.

In order to state Hain's results, let us first recall that a complex algebraic variety~$S$ is called a \emph{locally symmetric variety} if there is a semisimple algebraic group~$G$ over~$\bQ$, a maximal compact subgroup $K \subset G(\bR)$, and an arithmetic subgroup $\Gamma \subset G(\bQ)$, such that the symmetric space $G(\bR)/K$ is hermitian and $S = \Gamma\backslash G(\bR)/K$. As we have seen in~\ref{BasicPropSpecial}, special subvarieties are (at least after passage to a level cover) locally symmetric.

If $S_1 = \Gamma_1\backslash G_1(\bR)/K_1$ and $S_2 = \Gamma_2\backslash G_2(\bR)/K_2$ are locally symmetric varieties, a morphism $f \colon S_1 \to S_2$ is called a map of locally symmetric varieties if it is induced by a homomorphism of algebraic groups $G_1 \to G_2$ over~$\bQ$.

Hain calls a locally symmetric variety~$S$ \emph{good} if it has no locally symmetric divisors. Further, he mostly restricts his attention to locally symmetric varieties~$S$ for which the corresponding $\bQ$-group~$G$ is almost simple (i.e., $G^\ad$ is simple); for the problems that interest us this is no loss of generality. In case $G$ is almost simple, it can be shown that $S$ is not good only if the $\bQ$-rank of~$G$ is $\leq 2$ and if the non-compact factors of $G_\bR$ are all of the form $\SO(n,2)$ or $\SU(n,1)$, up to isogeny.

Define a \emph{locally symmetric family of abelian varieties} to be a principally polarized abelian scheme $X \to S$ such that $S$ is a locally symmetric variety and the corresponding morphism $S \to \sfA_g$ is a map of locally symmetric varieties. By a \emph{locally symmetric family of curves} we mean a curve $C\to S$ of compact type such that the corresponding relative Jacobian $J \to S$ is a locally symmetric family of abelian varieties.

\begin{theorem}[Hain]
\label{HainThm1}
Let $\pi\colon C \to S$ be a locally symmetric family of curves that is not isotrivial. Assume the $\bQ$-group that gives~$S$ is almost simple. Assume further that either $\pi$ is smooth, or $S$ is good and the generic fiber of~$\pi$ is smooth. Then $S$ is a quotient of a complex $n$-ball.
\end{theorem}

See \cite{Hain-Locsymm}, Theorem~1. 

It should be realized that for applications to Coleman's conjecture, the assumptions of the theorem are too strong. The main problem is that a special subvariety~$Z$ as in~\ref{ZinTorelli} only gives us (possibly after passing to a level cover) a locally symmetric family of abelian varieties $J \to Z$ such that the geometric generic fiber is a Jacobian. We may find a dominant morphism $Z^\prime \to Z$ such that the pullback $J^\prime \to Z^\prime$ is the relative Jacobian of a smooth curve $C \to Z^\prime$ but in general it is not possible to do this with $Z^\prime$ an open part of a locally symmetric variety. See \cite{Hain-Locsymm}, in particular Proposition~8.3.

Under weaker assumptions, Hain proves a second result.

\begin{theorem}[Hain]
\label{HainThm2}
Let $A \to S$ be a locally symmetric family of abelian varieties such that the morphism $p\colon S \to \sfA_g$ is not constant and factors through the Torelli locus~$\Torelli_g$. Write
\begin{align*}
S^\dec &:= \bigl\{s \in S(\bC) \bigm| \text{$A_s$ is the Jacobian of a singular curve,} \bigr\}\\
S^* &:= S(\bC) \setminus S^\dec\, ,\\
S^\he &:= \bigl\{s \in S(\bC) \bigm| \text{$A_s$ is the Jacobian of a hyperelliptic curve.} \bigr\}
\end{align*}
We assume the $\bQ$-group that gives~$S$ is almost simple, $S$ is good, and $S^* \neq \emptyset$. Then either $S$ is a ball quotient, or else $g\geq 3$, each component of~$S^\dec$ has complex codimension $\geq 2$ in~$S$, the family does not lift to a locally symmetric family of curves, and $S^* \cap S^\he$ is a non-empty divisor in~$S^*$, which moreover for $g > 3$ is nonsingular.
\end{theorem}

See \cite{Hain-Locsymm}, Theorem~2.

The proofs of Hain's results rely on a rigidity property of mapping class groups, which is a special case of the theorem of Farb and Masur in~\cite{Farb.Masur}.

De Jong and Zhang \cite{AJdJ-Zhang} have pushed Hain's results further, based on the observation that the hyperelliptic locus in~$\sfA_g\setminus \sfA_g^\dec$ is affine. 

\begin{theorem}[de Jong and Zhang]
\label{dJZThm}
With assumptions as in Theorem~\ref{HainThm2}, either $S$ is a ball quotient, or $S^\dec$ has codimension $\leq 2$ in~$S$, or the Baily-Borel compactification of~$S$ has a boundary of codimension~$\leq 2$.
\end{theorem}

\begin{corollary}
\label{dJZCor}
Let $S \subset \sfA_g$, with $g \geq 4$, be a Hecke translate of a Hilbert modular subvariety of~$\sfA_g$, i.e., $S$ is a special subvariety of PEL type obtained from a totally real field~$F$ of degree~$g$. Then $S$ is not contained in~$\Torelli_g$.
\end{corollary}

The Corollary was proved in \cite{AJdJ-Zhang}, except when $g=4$ and $F$ contains a quadratic subfield. That case was settled by Bainbridge and M\"oller in \cite{Bain-Moll}.

In addition to the results discussed here, there are several other results inspired by \ref{ZinTorelli}. Among these are papers of M\"oller, Viehweg and Zuo; see \cite{M.V.Z}, \cite{M.V.Z-Stability} and~\cite{V.Z}, and the recent paper~\cite{Andreatta-CO} by F.~Andreatta. These results support the Expectation~\ref{ZinTorelli}.

\mosubsection{Results of Ciliberto, van der Geer and Teixidor i Bigas}

In \cite{Cil.vdG.Teix} and \cite{Cil.vdG}, Ciliberto, van der Geer and Teixidor i Bigas have obtained some interesting results about the number of moduli of curves whose Jacobians have nontrivial endomorphisms. As we shall discuss below, when combined with the results of Hain and de Jong-Zhang, these results can be used to obtain some restrictions on the special subvarieties of PEL type that are contained in the Torelli locus.

As always in this section we work over~$\bC$.

\begin{theorem}[Ciliberto, van der Geer and Teixidor i Bigas, \cite{Cil.vdG.Teix}]
\label{CvdGTThm}
Let $Z \subset \sfM_g$, for $g \geq 2$, be an irreducible closed subvariety. Let $\Omega$ be an algebraic closure of the function field~$\bC(Z)$, let $C/\Omega$ be the curve corresponding to the geometric generic point of~$Z$, and assume that the Jacobian $J_C/\Omega$ has the property that $\bZ \subsetneq \End(J_C)$. Then $\dim(Z) \leq 2g-2$, and the intersection of~$Z$ with the hyperelliptic locus in~$\sfM_g$ has dimension at most~$g$.
\end{theorem}

Note that in this theorem it is not assumed that the image of~$Z$ in~$\sfA_g$ is a special subvariety.

In addition to the result as quoted here, the authors have some finer results about the case when $\dim(Z) \geq 2g-3$. For instance, they show that if $\dim(Z) = 2g-2$ then either $C$ is a cover of a non-constant elliptic curve $E$ over~$\Omega$, or $g=2$ and $\End^0(J_C)$ is a real quadratic field. (The case $\dim(Z) = 2g-3$ is analyzed by Ciliberto and van der Geer in~\cite{Cil.vdG}.)

\mosubsection{Excluding certain special subvarieties of PEL type}
One may use the results we have discussed to obtain restrictions on the special subvarieties $Z \subset \sfA_g$ of PEL type that can be contained in~$\Torelli_g$ with $Z \cap \Torelli_g^\open \neq \emptyset$. (We do not, however, see a way to apply such arguments to arbitrary special subvarieties.) Our result is as follows.

\begin{theorem}
Consider a a special subvariety $S \subset \sfA_g$ of PEL type, arising from PEL data $(D,*,V,\phi)$ as in \emph{Remark~\ref{PELData}}, with $\dim_\bQ(V) = 2g$.

{\upshape (i)} Suppose $D=F$ is a totally real field. (Albert Type~I.) If $g > 4$ then $S$ is not contained in~$\Torelli_g$.

{\upshape (ii)} Suppose $D$ is a quaternion algebra over a totally real field~$F$ that splits at all infinite places of~$F$. (Albert Type~II.) If $g >8$ then $S$ is not contained in~$\Torelli_g$.
\end{theorem}

For Albert's classification of the possible endomorphism algebras of abelian varieties we refer to \cite{MAV}, Section~21.

\begin{proof}
In both cases it follows from the results of \cite{Shimura-63} that $D$ equals the endomorphism algebra of the abelian variety corresponding to the geometric generic point of~$S$. Hence the generic abelian variety in our family is geometrically simple; in particular, if $S \subset \Torelli_g$ then $S \cap \Torelli_g^\open \neq \emptyset$.

First suppose $D=F$ is a totally real field. Let $e=[F:\bQ]$, which is an integer dividing~$g$. Then $S$ arises as a quotient of a product of $e$ copies of the Siegel space $\Siegel_h$ with $h = g/e$. This gives
\[
\dim(S) = e \cdot \frac{h(h+1)}{2} = \frac{g(g+e)}{2e}\, .
\]
Theorem~\ref{CvdGTThm} gives the inequality $g(g+e)/2e \leq 2g-2$, so we find
\[
e \geq \frac{g^2}{3g-4} > \frac{g}{3}\, .
\]
As $e$ divides~$g$, either $g$ is even and $e=g/2$ or $e=g$. However, for $g> 4$ the case $e=g$ is excluded by Corollary~\ref{dJZCor}. Assume then $g > 4$ is even and $e = g/2$; in this case $\dim(S) = 3g/2$. The assumptions in Hain's theorem~\ref{HainThm2} are satisfied and $S$ is not a ball quotient. Hence we obtain that the intersection with the hyperelliptic locus is a nonempty divisor and therefore has dimension $(3g/2)-1$. This contradicts Theorem~\ref{CvdGTThm}.

Next we consider the case where $D$ is a quaternion algebra over a totally real field~$F$ such that $D \otimes_\bQ \bR$ is isomorphic to a product of $e = [F:\bQ]$ copies of $M_2(\bR)$. In this case $2e$ divides~$g$ and $S$ is a quotient of a product of $e$ copies of~$\Siegel_h$ with $h = g/2e$; so $\dim(S) = e \cdot h(h+1)/2$ and $S$ is not a ball quotient. The boundary in the Baily-Borel compactification has codimension $g/2 > 2$. By Theorem~\ref{dJZThm} it therefore suffices to show that $S^\dec$ cannot have codimension $\leq 2$.

There are two possibilities for the geometric generic point of a component of~$S^\dec$. Either the corresponding abelian variety is isogenous to a product $X_1 \times X_2$, where $X_1$ and~$X_2$ both have an action by an order in~$D$. Straightforward calculation gives that in this case the codimension is at least $e\cdot (h-1)$. Using the relation $g=2eh$ we find that for $g>8$ the codimension is~$>2$. The other possibility is that there is a subfield $F^\prime \subset F$ and a quaternion algebra~$D^\prime$ with center~$F^\prime$ such that $D \cong F \otimes_{F^\prime} D^\prime$. With $e^\prime = [F^\prime:\bQ]$ the codimension in this case equals $(e-e^\prime) \cdot h(h+1)/2$. Straightforward checking of the possibilities shows that for $g>8$ this is greater than~$2$.
\end{proof}

It seems that similar arguments will also work for the Albert Type~III, and even for Type~IV we expect that we can obtain some non-trivial conclusions. We have not yet pursued this.

\mosubsection{The Schottky problem}
Expectation~\ref{ZinTorelli} is of course intimately related to the Schottky problem, which is the problem of characterizing which ppav $(A,\lambda)$ are Jacobians of curves. There are several solutions or conjectural solutions of this problem. We refer to the overview papers \cite{Beauv-Schottky}, \cite{Debarre-Schottky}, \cite{vdG-Schottky} for an introduction to this beautiful topic. It seems that none of the (conjectural) solutions discussed in these papers can be directly applied to problems such as Conjecture~\ref{ColemanConj} or Expectation~\ref{ZinTorelli}. Though for both sides---algebraic curves and their moduli on the one hand, abelian varieties and special subvarieties in the moduli space on the other hand---we have many techniques and results at our disposal, it is difficult to find a language, or a set of techniques, using which both sides can be described simultaneously. This is a difficulty that we believe lies at the heart of the matter.

\section{Examples of special subvarieties in the Torelli locus}
\label{ExaSsT}

In this section we discuss examples of special subvarieties $S \subset \Torelli_g$ of positive dimension, with $g \geq 4$ and $S \cap \Torelli_g^\open \neq \emptyset$. The examples we consider arise from families of cyclic covers of~$\bP^1$. Throughout this section we work over~$\bC$. 

We first look at a concrete example, following~\cite{AJdJ.Noot}. (The example was already given in~\cite{Shimura}.)

\begin{example}
\label{m=5a=1112}
Consider the family of curves given by
\[
y^5 = x(x-1)(x-t)\, ,
\]
where $t \in \bC\setminus\{0,1\}$ is a parameter. The complete and regular model~$C_t$ of this curve is a cyclic cover of~$\bP^1$ with group~$\mu_5$, with total ramification above $0$, $1$, $t$ and~$\infty$. The Hurwitz formula gives $\chi = 5 \cdot -2 + 4 \cdot (5-1) = 6$ so $g(C_t) = 4$. In the moduli space~$\sfA_4$ the corresponding family of Jacobians gives a $1$-dimensional subvariety, whose closure we call $Z \subset \sfA_4$. Clearly these Jacobians admit an action by $\bZ[\zeta_5]$, the ring of integers of the cyclotomic field $F := \bQ[\zeta_5]$. 

Consider the special subvariety $S \subset \sfA_4$ containing~$Z$ that is defined by this action of $\bZ[\zeta_5]$. More formally, fix a base point $b \in \bC\setminus\{0,1\}$ and choose a symplectic similitude $s\colon H_1(C_b,\bZ) \isomarrow V$ as in~\ref{AgShimVar}. Via~$s$, the Hodge structure on $H_1(C_b,\bQ)$ corresponds to a point $y \in \Siegel_4$ and we obtain the structure of an $F$-vector space on~$V_\bQ$. Let $M \subset G_\bQ = \CSp(V,\phi)_\bQ$ be the algebraic subgroup obtained as the intersection of $G$ with $\GL_F(V_\bQ)$. With notation as in Def.~\ref{SpecSubvarDef2}, $y$ lies in~$Y_M$. If we choose the base point~$b$ such that the Jacobian~$J_b$ is not of CM type (which is certainly possible) then there is a unique connected component $Y^+ \subset Y_M$ containing~$y$, and $S \subset \sfA_4$ is the special subvariety obtained as the image of this~$Y^+$ under the quotient map $\Siegel_4 \to \CSp(V,\phi)\backslash \Siegel_4 \isomarrow \sfA_4(\bC)$.

As we shall show, $\dim(S) = 1$. Assuming we know this, we conclude that $Z$, which is also $1$-dimensional and is contained in~$S$, contains an open dense subset of~$S$, in which case it follows from property~(b) in~\ref{BasicPropSpecial} that there are infinitely many values of~$t$ such that the Jacobian of~$C_t$ is of CM type. 

In order to calculate the dimension of~$S$, we need to know how $F$ acts on the tangent space of the Jacobian~$J_t$ at the origin. More precisely, $T_0(J_t)$ has the structure of a module over the ring $F \otimes_\bQ \bC = \prod_{\sigma\colon F \to\bC}\, \bC$. Hence we obtain a direct sum decomposition $T_0(J_t) = \oplus_\sigma\, T(\sigma)$. Let $n_\sigma$ denote the $\bC$-dimension of~$T(\sigma)$. These multiplicities~$n_\sigma$ do not depend on~$t$. By using the polarization, one can show that $n_\sigma + n_{\bar\sigma} = 2g/\varphi(5) = 2$ for all~$\sigma$, where $\bar\sigma$ denotes the complex conjugate of~$\sigma$. With this notation we have $H_\bR^\der \cong \prod \SU(n_\sigma,n_{\bar\sigma})$ and $\dim(S) = \sum n_\sigma n_{\bar\sigma}$, where the sum is taken over a set of representatives of the complex embeddings of~$F$ modulo complex conjugation. See~\cite{Shimura-63}, Theorem~5, and see below for further details.

As $T_0(J_t)$ is canonically dual to $H^0(C_t,\omega)$, we can calculate the multiplicities~$n_\sigma$ by writing down a basis of regular differentials on~$C_t$. (In fact, as we shall explain below the Chevalley-Weil formula gives a much quicker method.) In the example at hand, if $P_i$, for $i \in \{0,1,t,\infty\}$, is the unique point above~$i$, we have
\[
\Div(x) = 5P_0 - 5P_\infty\; ,\quad \Div(y) = P_0 + P_1 + P_t - 3P_\infty\; ,
\]
and
\[
\Div(dx) = 4P_0 + 4P_1 + 4P_t - 6P_\infty\, .
\]
As a basis of regular differentials we find
\[
\frac{dx}{y^2}\, ,\quad \frac{dx}{y^3}\, ,\quad \frac{dx}{y^4}\, ,\quad \frac{xdx}{y^4}\, .
\] 
The weights of the first two are dual, whereas the last two forms have the same weight. The conclusion, then, is that there is one pair $(\sigma,\bar\sigma)$ with multiplicities $(2,0)$, and one pair where the multiplicities are $(1,1)$. Hence $\dim(S) = 1\cdot 1 + 2\cdot 0 = 1$, and we conclude that $S = Z$ is a special subvariety contained in~$\Torelli_4$ with $S \cap \Torelli_4^\open \neq \emptyset$. In particular this proves that Coleman's conjecture~\ref{ColemanConj} does not hold for $g=4$.
\end{example}

\begin{remark}
\label{CMValuesUnknown}
In the above example, we do not know a method to decide for which values of the parameter~$t$ the Jacobian of~$C_t$ is a CM abelian variety. The same remark applies to the examples we shall discuss next. 
\end{remark}

\mosubsection{Families of cyclic covers}
\label{FamCycCov}
To obtain further such examples, the idea is to fix an integer $m \geq 2$, an integer $N \geq 4$, and monodromy elements $a_1,\ldots,a_N$ in $\bZ/m\bZ$; then we consider cyclic covers of~$\bP^1$ with group~$\mu_m$, branch points $t_1,\ldots,t_N$ in~$\bP^1$ and local monodromy $\exp(2\pi i a_j/m) \in \mu_m$ about~$t_j$. If the branch points are all in~$\bA^1$, this cover is given by the affine equation
\begin{equation}
\label{CycCovEq}
y^m = (x-t_1)^{a_1} (x-t_2)^{a_2} \cdots (x-t_N)^{a_N}\, ,
\end{equation}
with $\zeta \in \mu_m$ acting by $(x,y) \mapsto (x,\zeta \cdot y)$. Varying the branch points~$t_i$ gives us a family of curves, and it turns out that for certain choices of the data involved, the corresponding family of Jacobians traces out a special subvariety in~$\sfA_g$.

Write $a=(a_1,\ldots,a_N)$. The triple $(m,N,a)$ that serves as input for our construction has to satisfy some conditions. As already indicated, we want $m \geq 2$, as taking $m=1$ is clearly of no interest. Next, we require that $N \geq 4$. The reason is that from our family we obtain an $(N-3)$-dimensional subvariety in the moduli space (see below), and we are interested in special subvarieties of positive dimension. We assume that $a_i \not\equiv 0 \pmod{m}$ for all~$i$, and that $a_1 + \cdots  + a_N \equiv 0 \pmod{m}$. This means that the~$t_i$ are branch points and that there are no further branch points. Further we need to assume that the elements~$a_i$ generate the group $\bZ/m\bZ$ for otherwise the (smooth projective) curves given by~(\ref{CycCovEq}) are reducible.

We can find an open subscheme $T \subset (\bP^1)^N$, disjoint from the big diagonals, and a smooth proper curve $f \colon C \to T$ such that the fiber of~$f$ at a point $t = (t_1,\ldots,t_N)$ in $T(\bC)$ with all~$t_i$ in $\bA^1(\bC)$ is the complete regular curve given by~(\ref{CycCovEq}). The genus of these curves is given by
\begin{equation}
\label{gformula}
g = 1 + \frac{(N-2)m - \sum_{i=1}^N \gcd(a_i,m)}{2}\, .
\end{equation}
The relative Jacobian $J \to T$ then defines a morphism $\tau \colon T \to \sfA_g$, and we define $Z(m,N,a) \subset \sfA_g$ as the closure of the image of~$\tau$. Note that since we are only interested in the closure of~$\tau(T)$ there is no need to specify exactly which open subscheme $T \subset (\bA^1)^N$ we choose; see Remark~\ref{phiTSpecRem}.

We call two triples $(m,N,a)$ and $(m^\prime,N^\prime,a^\prime)$ equivalent if $m=m^\prime$ and $N=N^\prime$ and if the classes of $a$ and~$a^\prime$ in $(\bZ/m\bZ)^N$ are in the same orbit under $(\bZ/m\bZ)^* \times \gS_N$. Here $(\bZ/m\bZ)^*$ acts diagonally on $(\bZ/m\bZ)^N$ by multiplication, and the symmetric group~$\gS_N$ acts by permutation of the indices. The closed subvariety $Z(m,N,a) \subset \sfA_g$ only depends on the equivalence class of the triple $(m,N,a)$. 

The morphism $\tau \colon T \to \sfA_g$ factors through the quotient of~$T$ modulo the action of the group $\PGL_2(\bC) \times \gS_N$, with $\PGL_2(\bC) = \Aut(\bP^1)$ acting diagonally on $(\bP^1)^N$ and $\gS_N$ acting by permutation of the diagonals. (Without loss of generality we may assume $T \subset (\bP^1)^N$ is stable under $\PGL_2(\bC) \times \gS_N$.) The subvariety $Z(m,N,a) \subset \sfA_g$ has dimension $N-3$. Note that we may also fix three of the branch points to lie at $0$, $1$ and~$\infty$, as we did in Example~\ref{m=5a=1112}; this has the effect of replacing~$T$ with a closed subvariety of dimension $N-3$ on which the morphism~$\phi$ is generically finite. For instance, the example considered in~\ref{m=5a=1112} corresponds to the triple $(m,N,a) =  \bigl(5,4,(1,1,1,2)\bigr)$.

The Jacobians $J_t$ in our family come equipped with an action of the group ring $\bZ[\mu_m]$. Let $S(\mu_m) \subset \sfA_g$ be the special subvariety containing $Z(m,N,a)$ that is defined by this action. (In order to make this precise we follow the recipe given in Example~\ref{m=5a=1112}.) In all cases we have the inequality
\begin{equation}
\label{N-3<dim}
N-3 \leq \dim S(\mu_m)\, ,
\end{equation}
and if we have $N-3 = \dim S(\mu_m)$ then the conclusion is that $Z(m,N,a) = S(\mu_m)$ is a special subvariety in the Torelli locus that meets the open Torelli locus. (If $N-3 < \dim S(\mu_m)$ then a priori we know nothing; see \ref{ExcludeExa} below.) The question is therefore how to calculate the dimension of $S(\mu_m)$. Before we discuss this in general, let us look at another example.

\begin{example}
\label{m=9a=1116}
Consider the family of curves given by
\[
y^9 = x(x-1)(x-t)\, .
\]
With notation as above we have $m=9$ and $N=4$, with local monodromy about the branch points given by $a = (1,1,1,6)$.

The complete regular model~$C_t$ is a cyclic cover of~$\bP^1$ with group~$\mu_9$. The points $0$, $1$ and~$t$ are totally ramified, so we have unique points~$P_0$, $P_1$, and~$P_t$ above them. There are three points $P_\infty^{(1)}$, $P_\infty^{(2)}$ and~$P_\infty^{(3)}$ above~$\infty$, each with ramification index~$3$. The Hurwitz formula gives $\chi = 9 \cdot -2 + 3 \cdot 8 + 3 \cdot 2 = 12$ so $g=7$, which agrees with~(\ref{gformula}).

The Jacobian~$J_t$ contains as an isogeny factor the Jacobian of the curve~$C^\prime_t$ given by $u^3 = x(x-1)(x-t)$. By a similar calculation, $C^\prime_t$ has genus~$1$. Its Jacobian~$J_t^\prime$ is an elliptic curve with complex multiplication by $\bZ[\zeta_3]$, and the family of Jacobians~$J_t^\prime$ is isotrivial over $\bP^1\setminus\{0,1\}$. (Note that in the given equation for~$C^\prime_t$ there is no ramification over~$\infty$, so $C^\prime_t$ is geometrically isomorphic to the curve given by $s^3 = v(v-1)$. More explicitly, let $\tilde{T} \to \bP^1\setminus\{0,1\}$ be the cyclic cover of degree~$3$ given by the equation $\mu^3 = t(t-1)$. After base change to~$\tilde{T}$ the family of curves $C^\prime_t$ becomes isomorphic to the constant curve defined by $s^3 = v(v-1)$ via the isomorphism given by $s = \mu u/(x-t)$ and $v = (1-t)x/(x-t)$. )

Let $J_t^\new$ (the ``new part'') be the quotient of $J_t$ modulo~$J_t^\prime$; this is an abelian variety of dimension~$6$ on which we have an action by  $\bZ[\zeta_9]$. In order to determine the dimension of the special subvariety in~$\sfA_6$ given by this action, we again calculate the multiplicities of the action on the tangent space. In this case, we already know in advance that there will be one non-primitive character of~$\mu_9$ occurring in $H^0(C_t,\omega)$. We have
\[
\Div(x) = 9P_0 - 3P^{(1)}_\infty - 3P^{(2)}_\infty - 3P^{(2)}_\infty\; ,\quad \Div(y) = P_0 + P_1 + P_t - P^{(1)}_\infty - P^{(2)}_\infty - P^{(3)}_\infty\; ,
\]
and
\[
\Div(dx) = 8P_0 + 8P_1 + 8P_t - 4P^{(1)}_\infty - 4P^{(2)}_\infty - 4P^{(3)}_\infty\, .
\]
As a basis of regular differentials we find
\[
\frac{dx}{y^4}\, ,\quad \frac{dx}{y^5}\, ,\quad \frac{dx}{y^6}\, ,\quad \frac{dx}{y^7}\, ,\quad \frac{dx}{y^8}\, ,\quad \frac{xdx}{y^7}\, ,\quad \frac{xdx}{y^8}\, .
\]
As predicted, this gives one non-primitive character of~$\mu_9$ (corresponding to $dx/y^6$); for the rest we find one complex conjugate pair with multiplicities $(1,1)$, and two complex conjugate pairs for which the multiplicities are $(2,0)$. The special subvariety in~$\sfA_6$ given by the action of $\bQ[\zeta_9]$ with this collection of multiplicities has dimension $1\cdot 1 + 2\cdot 0 + 2\cdot 0 = 1$. It follows that the $J_t^\new$ trace out a dense open subset of this special subvariety. As the original Jacobians~$J_t$ are isogenous to $J_t^\prime \times J_t^\new$ with $J_t^\prime$ an isotrivial family with fibers of CM type, we conclude that the family of~$J_t$ traces out in~$\sfA_7$ a (dense subset of a) special subvariety. As before, this implies that there are infinitely many values of~$t$ for which $J_t$ is of CM type. (The $t$ for which this happens are even analytically dense in~$\bP^1(\bC)$.) 
\end{example}

\begin{remark}
\label{DimShimCalc}
Let $(M,Y)$ be a Shimura datum. The dimension of the associated Shimura varieties $\Shim_K(M,Y)$ only depends on the structure of the real adjoint group~$M_\bR^\ad$. Indeed, the dimension of the Shimura variety equals the complex dimension of the space~$Y$, and the components of~$Y$ can be described as hermitian symmetric domains associated with the connected Lie group $M^\ad(\bR)^+$, where the superscript ``$+$'' denotes the identity component for the analytic topology. In particular, if $M^\ad_\bR \cong Q_1 \times \cdots \times Q_l$ is the decomposition of~$M^\ad_\bR$ as a product of $\bR$-simple groups, the dimension of $\Shim_K(M,Y)$ can be calculated as a sum $\delta(Q_1) + \cdots + \delta(Q_l)$, where the contribution $\delta(Q_i)$ of the factor~$Q_i$ can be looked up in \cite{Helgason}, Table~V. (Caution: the dimensions given there are the real dimensions.) The cases most relevant for our discussion are the following.
\begin{enumerate}
\item $\delta(Q_i) = 0$ if $Q_i$ is anisotropic, i.e., if $Q_i(\bR)$ is compact;
\item $\delta(Q_i) = h(h+1)/2$ if $Q_i \cong \PSp_{2h,\bR}$;
\item $\delta(Q_i) = pq$ if $Q_i \cong \PSU(p,q)$.
\end{enumerate}
\end{remark}

\mosubsection{Calculating the dimension of $S(\mu_m)$}
We return to the general situation of a family of curves $f\colon C \to T$ associated with the data $(m,N,a)$. Our next goal is to calculate the dimension of the special subvariety $S(\mu_m) \subset \sfA_g$ that contains $Z = Z(m,N,a)$. Choose a Hodge-generic base point $b \in T(\bC)$, choose a symplectic similitude $s \colon H_1(C_b,\bZ) \isomarrow V$ as in~\ref{AgShimVar}, and let $y \in \Siegel_g$ be the point corresponding to the Hodge structure on $H_1(C_b,\bZ)$, viewed as a Hodge structure on~$V$ via~$s$.

The Jacobian $J \to T$ comes equipped with an action of the group ring $\bZ[\mu_m]$. The Hodge classes we want to have over~$S(\mu_m)$ are these endomorphisms, which means we are in the situation of Example~\ref{SpecSubvarEndom}. The algebraic group $M \subset \CSp(V_\bQ,\phi)$ we need to consider is given by
\[
M := \CSp(V_\bQ,\phi) \cap \GL_{\bQ[\mu_m]}(V_\bQ)\, .
\]

We can  calculate the dimension of $S(\mu_m)$  via a deformation argument. In other words, if $Y_M \subset \Siegel_g$ is as defined just before Definition~\ref{SpecSubvarDef2}, we calculate the dimension of the tangent space of~$Y_M$ at te point~$y$. 

The Hodge structure on~$V_\bQ$ is of type $(-1,0) + (0,-1)$, and the polarization gives us a symplectic form $\phi \colon V_\bQ \times V_\bQ \to \bQ(1)$. As discussed in Remark~\ref{BorelEmb}, the Hodge filtration $\Fil^0 = V_\bC^{0,-1} \subset V_\bC$ is a Lagrangian subspace, so we have an induced isomorphism $\bar{\phi} \colon \Fil^0 \isomarrow V_\bC/\Fil^0$. The tangent space of~$\Siegel_g$ at the point~$y$ is given by~\eqref{TySiegel}.

The extra structure we now have is an action of the algebra $D := \bQ[\mu_m]$ on~$V_\bQ$. Let $d\mapsto \bar{d}$ denote the involution of~$D$ induced by the inversion in the group~$\mu_m$. The form~$\phi$ has the property that $\phi(dv,v^\prime) = \phi(v,\bar{d}v^\prime)$ for all $d \in D$ and $v$, $v^\prime \in V_\bQ$. 

The $D$-action on~$V_\bQ$ induces on $V_\bC^{-1,0}$ and~$V_\bC^{0,-1}$ the structure of a module over the ring
\[
D \otimes_\bQ \bC = \prod_{n \in \bZ/m\bZ}\, \bC\, .
\]
Correspondingly, we have decompositions
\[
V_\bC/\Fil^0 = V_\bC^{-1,0} = \bigoplus_{n \in \bZ/m\bZ}\, V_{\bC,(n)}^{-1,0}
\quad\mbox{and}\quad
\Fil^0 = V_\bC^{0,-1} = \bigoplus_{n \in \bZ/m\bZ}\, V_{\bC,(n)}^{0,-1}\, .
\]
The involution of $D_\bC$ obtained by linear extension of the involution $d \mapsto \bar{d}$ exchanges the factors~$\bC$ indexed by the classes $n$ and~$-n$. It follows that the perfect pairing $\bar{\phi} \colon \Fil^0 \times (V_\bC/\Fil^0) \to \bC$ restricts to perfect pairings
\[
\bar{\phi}_n \colon V_{\bC,(n)}^{0,-1} \times V_{\bC,(-n)}^{-1,0} \to \bC\, .
\]

For $n \in \bZ/m\bZ$, define
\[
d_n := \dim_\bC V_{\bC,(n)}^{0,-1}\, .
\]
Note that $d_{(0 \bmod m)} = 0$, i.e., $V_{\bC,(0\bmod m)} = (0)$, as the $\mu_m$-invariant subspace in~$V$ is the first homology of the base curve~$\bP^1$, which is zero. We shall see in Proposition~\ref{ChevWeil} below how to calculate the~$d_n$ in terms of the given data $(m,N,a)$.

The tangent space $T_y(Y_M) \subset T_y(\Siegel_g)$ is the $\bC$-subspace of elements $\beta \in \Hom^\sym(\Fil^0,V_\bC/\Fil^0)$ that are $D_\bC$-linear. Any such~$\beta$ can be written as $\beta = \sum\, \beta_n$, where the $\beta_n \colon V_{\bC,(n)}^{0,-1} \to V_{\bC,(n)}^{-1,0}$ are $\bC$-linear maps that satisfy
\begin{equation}
\label{betanRel}
\bar{\phi}_n\bigl(v,\beta_{-n}(v^\prime)\bigr) = \bar{\phi}_{-n}\bigl(v^\prime,\beta_n(v)\bigr)\quad \mbox{for all $v \in V_{\bC,(n)}^{0,-1}$ and $v^\prime \in V_{\bC,(-n)}^{0,-1}$.}
\end{equation}
If $n \not\equiv -n \pmod{m}$, this last condition gives a duality between $\beta_{-n}$ and~$\beta_n$; this means the linear map~$\beta_n$ can be chosen arbitrarily and $\beta_{-n}$ is determined by~$\beta_n$. The situation is different if $n \equiv -n \pmod{m}$. Of course, this only occurs (with $n \not\equiv 0$) if $m=2k$ is even. In this case we have a perfect pairing $\bar{\phi}_k \colon V_{\bC,(k)}^{0,-1} \times V_{\bC,(k)}^{-1,0} \to \bC$, and \eqref{betanRel} gives
\begin{equation*}
\begin{split}
\beta_k \in {} & \Hom^\sym\bigl(V_{\bC,(k)}^{0,-1},V_{\bC,(k)}^{-1,0}\bigr) \\ 
& := \bigl\{\beta \colon V_{\bC,(k)}^{0,-1} \to V_{\bC,(k)}^{-1,0} \bigm| \bar{\phi}_k(v,\beta_k(v^\prime)) = \bar{\phi}_k(v^\prime,\beta_k(v))\ \mbox{for all $v$, $v^\prime \in V_{\bC,(k)}^{0,-1}$} \bigr\}\, .
\end{split}
\end{equation*}
(Cf.\ \eqref{TySiegel}.) We find that the dimension of $T_y(Y_M)$ equals
\[
\sum_{\genfrac{}{}{0pt}{1}{\pm n \in (\bZ/m\bZ)/\{\pm 1\}}{2n \not\equiv 0 \pmod{m}}}\!\! d_{-n} d_n\; + \begin{cases} \frac{d_k \cdot \bigl(d_k +1\bigr)}{2} & \text{if $m=2k$ is even;}\\ 0 & \text{if $m$ is odd.}\end{cases}
\]

As the final step in the calculation we need to calculate the dimensions~$d_n$ in terms of the given triple $(m,N,a)$. The result is as follows.

\begin{proposition}[Hurwitz, Chevalley-Weil]
\label{ChevWeil}
The dimensions $d_n := \dim_\bC V_{\bC,(n)}^{0,-1}$ are given by $d_n = 0$ if $n \equiv 0 \pmod{m}$ and
\begin{equation}
\label{ChevWeilForm}
d_n = -1 + \sum_{i=1}^N \fracpart{\frac{-na_i}{m}} 
\quad \mbox{if $n \not\equiv 0 \pmod{m}$,}
\end{equation}
where $\fracpart{x} = x - \lfloor x\rfloor$ denotes the fractional part of a number~$x$.
\end{proposition}

We note that $V_\bC^{0,-1}$ is naturally isomorphic to the space $H^0(C_b,\Omega^1)$ of global differentials on the curve. The given formula is then a special case of a classical result by Chevalley and Weil about the structure of $H^0(C_b,\Omega^1)$ as a representation of the Galois group~$\mu_m$, which in the cyclic case is already due to Hurwitz. See \cite{Morr.Pinkh} Section~3 for a modern proof. Another proof, using the holomorphic Lefschetz formula, can be found in \cite{Schoen}, Lemma~1.6b.

Putting everything together, the dimension of $S(\mu_m)$ is given by the following result.

\begin{proposition}
\label{DimShim}
Consider a triple $(m,N,a)$ as in Section~\ref{FamCycCov}. Then the dimension of the special subvariety $S(\mu_m)$ that contains $Z(m,N,a)$ is given by
\[
\dim S(\mu_m) = \sum\, d_{-n} d_n\; + \begin{cases} \frac{d_k \cdot \bigl(d_k +1\bigr)}{2} & \text{if $m=2k$ is even;}\\ 0 & \text{if $m$ is odd,}\end{cases}
\]
where the sum runs over the pairs $\pm n$ in $\bZ/m\bZ$ with $2n \not\equiv 0 \pmod{m}$, and where the $d_n$ are given by~\eqref{ChevWeilForm}.
\end{proposition}

\begin{remark}
Instead of using a tangent space computation, we may calculate the dimension of~$S(\mu_m)$ by analyzing the real adjoint group~$M_\bR$; see Remark~\ref{DimShimCalc} and~\cite{Moonen-CycCov}. It can be shown that
\[
M_\bR^\ad \cong \prod \PSU(d_n,d_{-n}) \times \begin{cases} \PSp_{2d_k,\bR} & \text{if $m=2k$ is even;}\\ \{1\} & \text{if $m$ is odd,}\end{cases}
\]
where the first product runs over the pairs $\pm n$ in $\bZ/m\bZ$ with $2n \not\equiv 0 \pmod{m}$.
\end{remark}

\mosubsection{An inventory of known special subvarieties in the Torelli locus}
Now that we have an explicit formula for the dimension of $S(\mu_m)$ in terms of the given data, it is easy to check, in each given example, if in the inequality~\eqref{N-3<dim} we have an equality. Recall that if $N-3 = \dim S(\mu_m)$, we conclude that $Z(m,N,a)$ is dense in the special subvariety $S(\mu_m) \subset \sfA_g$, in which case it follows that the family of Jacobians $J \to T$ (as in section~\ref{FamCycCov}) contains infinitely many different Jacobians of CM type.

In order to illustrate how well this works, let us redo Example~\ref{m=5a=1112}. We have $(m,N,a) = \bigl(5,4,(1,1,1,2)\bigr)$. The Hurwitz-Chevalley-Weil formula gives
\[
d_1 = -1 + 3 \cdot \fracpart{\frac{-1}{5}} + \fracpart{\frac{-2}{5}} = 2\, ,
\]
and in a similar way we find 
\[
d_2 = 1\, ,\quad d_3 = 1\, ,\quad d_4 = 0\, ,
\]
which agrees with the basis of regular differentials we have found. Proposition~\ref{DimShim} gives $\dim S(\mu_m) = 1 = N-3$, and as before we conclude that $Z(m,N,a)$ is a special subvariety.

Using a computer it is not hard to do a systematic search for triples $(m,N,a)$ with $N-3 = \dim S(\mu_m)$. Up to equivalence (as defined in Section~\ref{FamCycCov}), one finds twenty such triples. They are listed in Table~\ref{TableExamples}. Moreover, it was proven by Rohde in~\cite{Rohde} that these are the only triples, up to equivalence, for which $N-3=\dim S(\mu_m)$.

\begin{table}[ht]
\setlength\extrarowheight{4pt}
\caption{Examples of special subvarieties in the Torelli locus}
\label{TableExamples}
\begin{minipage}{.5\textwidth}
\begin{center}
\begin{tabular}{|c|c|c|c|c|}
\hline
& genus & $m$ & $N$ & $a$ \\
\hline
\nr & 1 & 2 & 4 & (1,1,1,1) \cr
\hline
\nr & 2 & 2 & 6 & (1,1,1,1,1,1) \cr
\hline
\nr & 2 & 3 & 4 & (1,1,2,2) \cr
\hline
\nr & 2 & 4 & 4 & (1,2,2,3) \cr
\hline
\nr & 2 & 6 & 4 & (2,3,3,4) \cr
\hline
\nr & 3 & 3 & 5 & (1,1,1,1,2) \cr
\hline
\nr & 3 & 4 & 4 & (1,1,1,1) \cr
\hline
\nr & 3 & 4 & 5 & (1,1,2,2,2) \cr
\hline
\nr & 3 & 6 & 4 & (1,3,4,4) \cr
\hline
\nr & 4 & 3 & 6 & (1,1,1,1,1,1) \cr
\hline
\end{tabular}
\end{center}
\end{minipage}
\begin{minipage}{.49\textwidth}
\begin{center}
\begin{tabular}{|c|c|c|c|c|}
\hline
& genus & $m$ & $N$ & $a$ \\
\hline
\nr & 4 & 5 & 4 & (1,3,3,3) \cr
\hline
\nr & 4 & 6 & 4 & (1,1,1,3) \cr
\hline
\nr & 4 & 6 & 4 & (1,1,2,2) \cr
\hline
\nr & 4 & 6 & 5 & (2,2,2,3,3) \cr
\hline
\nr & 5 & 8 & 4 & (2,4,5,5) \cr
\hline
\nr & 6 & 5 & 5 & (2,2,2,2,2) \cr
\hline
\nr & 6 & 7 & 4 & (2,4,4,4) \cr
\hline
\nr & 6 & 10 & 4 & (3,5,6,6) \cr
\hline
\nr & 7 & 9 & 4 & (3,5,5,5) \cr
\hline
\nr & 7 & 12 & 4 & (4,6,7,7) \cr
\hline
\end{tabular}
\end{center}
\end{minipage}
\end{table}

\mosubsection{Excluding further examples}
\label{ExcludeExa}
In the discussion so far, the argument is based on the fact that $Z(m,N,a) \subset \sfA_g$ is visibly contained in the special subvariety $S(\mu_m)$, and we look for examples where the two have the same dimension. If in some given case we find that $\dim Z(m,N,a) = N-3 < \dim S(\mu_m)$, this does not, a priori, imply that $Z(m,N,a)$ is not special. Put differently, in addition to the endomorphisms in $\bZ[\mu_m]$, there might be Hodge classes in our family of Jacobians $J \to T$ that we just happen not to see. Even in individual examples, it is usually not so easy to exclude this. 

One method to do this was given by de Jong and Noot in~\cite{AJdJ.Noot}, Section~5; they prove there that for $m > 7$ not divisible by~$3$ the family of curves given by $y^m = x(x-1)(x-t)$, which in our language corresponds to the triple $(m,N,a) = \bigl(m,4,(1,1,1,m-3)\bigr)$, does \emph{not} give a special subvariety. The method is based on results of Dwork and Ogus in~\cite{Dwork.Ogus}.

Extending this to arbitrary families, it was proven in~\cite{Moonen-CycCov} that the twenty examples we have found are the only ones such that the image $Z(m,N,a) \subset \sfA_g$ is a special subvariety.

\begin{theorem}
\label{BM20Fams}
Consider data $(m,N,a)$ as in~\emph{\ref{FamCycCov}}. Then $Z(m,N,a) \subset \sfA_g$ is a special subvariety if and only if $(m,N,a)$ is equivalent to one of the twenty triples listed in \emph{Table~\ref{TableExamples}}.
\end{theorem}

\begin{remark}
Examples (6), (8), (10), (11), (16) and~(17) in Table~\ref{TableExamples} were given by Shimura in~\cite{Shimura}. Examples (10), (11) and~(17) were given, in a language that is closer to the present paper, by de Jong and Noot in~\cite{AJdJ.Noot}; they also explained the relevance of such examples for Coleman's conjecture. (The fact that these examples already occurred in~\cite{Shimura} was recognized only later. In connection with this, note that \cite{Shimura} appeared more than twenty years before Coleman stated his conjecture.)

In retrospect, it is surprising that the complete list of examples was obtained only recently. Example~(19) was found by one of us in 2003; see~\cite{FO-ICSNote}. All twenty examples were found by Rohde in~\cite{Rohde} (Example~(2) is somewhat hidden there, but see op.\ cit.\ Corollary 5.5.2) and independently by one of us, with the help of a small computer program. See~\ref{FamsCoverGeneral} below for some further examples, coming from families of covers with a non-cyclic Galois group.

It should be mentioned that there is some connection between these examples and the theory of Deligne and Mostow in \cite{Del.Mostov}, \cite{Mostov1988}; see also Looijenga's overview paper~\cite{EL-Unif}. The relation is that, in the setting of~\ref{FamCycCov}, the family of Jacobians $J \to T$ has a ``new part'' (cf.\ Example~\ref{m=9a=1116}), and that the monodromy of the corresponding VHS can be described as the monodromy on a space of hypergeometric functions. In this context there is an easy criterion for the arithmeticity of the monodromy group; see \cite{Del.Mostov}, Proposition~12.7 or \cite{EL-Unif}, Theorem~4.3. As far as we know there is, however, no easy way to obtain from this, and the resulting tables in \cite{Del.Mostov} and~\cite{Mostov1988}, a classification result such as in Theorem~\ref{BM20Fams}. Apart from the fact that the Deligne-Mostow arithmeticity criterion only applies under some condition on the local monodromy elements (condition (INT) of \cite{Del.Mostov}, p.~25) that in our situation is not always satisfied, it only gives us information about certain direct summands of the VHS that we want to study, and as already explained in Remark~\ref{A1xA2Rem} we in general need more.
\end{remark}

\section{Some questions}

In this section, with the exception of \ref{cancoord} and~\ref{analogy}, the base field is~$\bC$. Most questions below could also be formulated over an algebraic closure of~$\bQ$.

\begin{question} 
\emph{How do we construct curves for which the Jacobian is a CM abelian variety?} If we have a special subvariety in~$\Torelli_g$ that intersects~$\Torelli_g^\open$, the existence of CM Jacobians follows from property~(b) in~\ref{BasicPropSpecial}. Note that in such cases we typically have no control over which fibers in our family are the CM fibers; cf.\ Remark~\ref{CMValuesUnknown}.

On the other hand, we could look for curves that have many automorphisms. For instance, the Fermat curves $X^n+Y^n+Z^n=0$ in~$\bP^2$ have Jacobians of CM type. Similarly, cyclic covers of~$\bP^1$ with $3$ branch points have CM Jacobians. For a given $g \geq 2$ there are only finitely many such covers of genus~$g$, however. Further such examples are given in \cite{FO-ModCM},~5.15.

In \cite{FO-ModCM}, a curve~$C$ is called a curve \emph{with many automorphisms} if the deformation functor of the pair $\bigl(C,\Aut(C)\bigr)$ is a $0$-dimensional scheme. In many cases such a curve has a CM Jacobian.

These are the only methods known to us to construct, or to prove the existence of, CM Jacobians.
\end{question}

\begin{question}
\emph{Do we know the existence of, or can we construct, a curve~$C$ of genus at least~$4$ such that $\Aut(C) = \{\id\}$ and such that the Jacobian of~$C$ is an abelian variety of CM type?}
\end{question}

\begin{question}
\emph{Does there exist $g>3$ and a special subvariety $Z \subset \cA_g$ contained in the Torelli locus~$\Torelli_g$ such that the geometric generic fiber over $T$ gives an abelian variety with endomorphism ring equal to $\bZ$~?} We do not know a single example.
\end{question}

\begin{remark}
In \ref{ZinTorelli} we have seen an expectation about the (non-)existence of special subvarieties in the Torelli locus. In order to make the question more precise, consider, for $g \in \bZ_{>0}$ the set 
\[
\SpecialTor(g) := \bigl\{\text{special subvarieties $Z \subset \Torelli_g$ with $\dim(Z) > 0$ and $Z \cap \Torelli_g^\open \neq \emptyset$} \bigr\}
\]
of special subvarieties of positive dimension, contained in the Torelli locus, and not fully contained in the boundary of~$\Torelli_g$. The expectation is that for $g \gg 0$ we have $\SpecialTor(g) = \emptyset$; see \ref{ZinTorelli}.

We would like to classify all pairs $(g,Z)$  with  $Z \in  \SpecialTor(g)$.  For $g=2$ and $g=3$ we have $\Torelli_g =  \sfA_g$ and in this case every special subvariety of $\sfA_{g}$ is of PEL type; see~\cite{BM.YuZ-MathAnn}. Hence in this case we can classify all pairs $(g,Z)$, up to Hecke translation, by listing all possible endomorphism algebras.  We know that $\SpecialTor(g) \neq \emptyset$ for all $g < 8$. However, already for $g=4$ we do not have a good description of $\SpecialTor(4)$. It seems very difficult to describe $\SpecialTor(g)$ for arbitrary~$g$.
\end{remark}

\mosubsection{Non-PEL Shimura curves for $g=4$}
\label{Mumf4} 
In \cite{Mumf-ShimPaper}, \S~4, Mumford shows there exist $1$-dimensional special subvarieties $Z \subset \sfA_4$ that are not of PEL type. The abelian variety corresponding to the geometric generic fiber of~$Z$ has endomorphism algebra~$\bZ$. The curves~$Z$ are complete. Note that $\Torelli_4 \subset \sfA_4$ is a closed subvariety of codimension one, which by a result of Igusa \cite{Igusa} is ample as a divisor. Hence we see that $Z \cap \Torelli_4 \neq \emptyset$. 

\begin{question}
\label{MumfCurves}
\emph{Is there a ``Mumford curve'' $Z \subset \sfA_4$ that is contained in $\Torelli_4$~?}
\end{question}

If $Z \subset \sfA_4$ is a $1$-dimensional special subvariety of the type constructed by Mumford, the abelian variety corresponding to the geometric generic point of~$Z$ has endomorphism ring~$\bZ$; hence if $Z \subset \Torelli_4$ then $Z$ meets~$\Torelli_4^\open$.

The examples constructed by Mumford, and generalizations thereof, have been studied in detail by Noot in \cite{Noot-ellGal} and~\cite{Noot-MumfFam}. In particular, \cite{Noot-MumfFam}, Section~3 contains a detailed analysis of the possible CM points on the special curves $Z \subset \sfA_4$ constructed by Mumford. In particular, it is shown there that the (geometric) CM fibers are either absolutely simple or are isogenous to a product $E \times Y$ with $E$ an elliptic curve, $Y$ an abelian threefold, and $\End^0(E)$ isomorphic to a subfield of $\End^0(Y)$. It might be possible to use this to get some non-trivial information related to the above question. In connection with what was discussed in~\ref{ModifColeman}, let us note that there are no Weyl type CM fibers in these families.

\begin{question}
\emph{For which $g \geq 2$ does there exist a positive dimensional subvariety $Z\subset \Torelli_g$ with $Z \cap \Torelli_g^\open \neq \emptyset$, such that the abelian variety corresponding with the geometric generic point of~$Z$ is isogenous to a product of elliptic curves?}  

This question was stimulated by results in \cite{Kukulies},  which say that under more restrictive conditions such a family does not exist for large~$g$.
\end{question}

\begin{question}
\label{FamsCoverGeneral}
In Section~\ref{ExaSsT} we have seen examples of special subvarieties $Z \subset \Torelli_g$ with $Z \cap \Torelli_g^\open \neq \emptyset$ arising from families of cyclic covers of~$\bP^1$. \emph{Can one obtain further such examples by taking non-cyclic covers, or from a family of covers of another base curve?}

To make this more precise, consider a (complete, nonsingular) curve~$B$ over~$\bC$ of genus~$h$, and let $\cG$ be a finite group. Define a group $\Pi = \Pi(h,N)$ by
\[
\Pi := \bigl\langle \alpha_1,\ldots,\alpha_h,\beta_1,\ldots,\beta_h,\gamma_1,\ldots,\gamma_N\bigm| [\alpha_1,\beta_1] \cdots [\alpha_h,\beta_h] \cdot \gamma_1 \cdots \gamma_N = 1\bigr\rangle\, .
\]
Given $t_1,\ldots,t_N$ in~$B$, fix a presentation 
\[
\pi_1\bigl(B\setminus\{t_1,\ldots,t_N\}\bigr) \isomarrow \Pi\, ,
\]
and fix a surjective homomorphism $\psi \colon \Pi \twoheadrightarrow \cG$ such that $\psi(\gamma_i) \neq 1$ for all $i=1,\ldots,N$. Correspondingly, we have a Galois cover $C_t \to B$ with group~$\cG$, branch points $t_1,\ldots,t_N$ in~$B$, and with local monodromy about~$t_i$ given by the element~$\psi(\gamma_i)$. Varying the branch points, we get a family of curves $C \to T$, for some open $T \subset B^N$. The corresponding family of Jacobians gives a moduli map $T \to \sfA_g$; denote the image by $Z^\open$, with Zariski closure $Z \subset \Torelli_g \subset \sfA_g$. Having set the scene in this way we can ask for which choices of the data involved $Z$ is a special subvariety of positive dimension.
\begin{itemize}
\item[(a)] \emph{Are there examples with non-cyclic Galois group such that $Z \subset \sfA_g$ is a special subvariety of positive dimension?}
\item[(b)] \emph{Are there examples where $B$ is not a rational curve, and such that $Z \subset \sfA_g$ is a special subvariety of positive dimension?} Note that if $Z$ is special, the Jacobian of~$B$ is a CM abelian variety.
\end{itemize}

As for question~(a), we do have some examples that are obtained from families of covers of~$\bP^1$ with a non-cyclic abelian Galois groups. The examples presently known to us are listed in Table~\ref{TableExaNonCyc}. In these examples we consider families of covers of~$\bP^1$ with Galois group of the form $A = (\bZ/m_1\bZ) \times (\bZ/m_2\bZ)$ with $m_1|m_2$, with $N \geq 4$ branch points, and with local monodromy about the branch points~$t_i$ given by an $N$-tuple $a = (a_1,\ldots,a_N)$ in~$A^N$. This gives rise to an $(N-3)$-dimensional closed irreducible subvariety $Z(m,N,a)$, where now $m=(m_1,m_2)$. For the Jacobians~$J_t$ in our family we have a (generally non-injective) homomorphism $\bQ[A] \to \End^0(J_t)$; this defines a special subvariety $S(m) \subset \sfA_g$ of PEL type with $Z(m,N,a) \subseteq S(m)$. We calculate $\dim\bigl(S(m)\bigr)$, and if this dimension equals $N-3$ then we conclude that $Z(m,N,a)$ is special.

\begin{table}[ht]
\setlength\extrarowheight{4pt}
\caption{Examples of special subvarieties in the Torelli locus (continued)}
\label{TableExaNonCyc}
\begin{center}
\begin{tabular}{|c|c|c|c|c|}
\hline
& genus & group & $N$ & $a$ \\
\hline
\nr & 1 & $(\bZ/2\bZ) \times (\bZ/2\bZ)$ & 4 & $\bigl((1,0),(1,0),(0,1),(0,1)\bigr)$ \cr
\hline
\nr & 3 & $(\bZ/2\bZ) \times (\bZ/4\bZ)$ & 4 & $\bigl((1,0),(1,1),(0,1),(0,2)\bigr)$ \cr
\hline
\nr & 3 & $(\bZ/2\bZ) \times (\bZ/4\bZ)$ & 4 & $\bigl((1,0),(1,2),(0,1),(0,1)\bigr)$ \cr
\hline
\nr & 4 & $(\bZ/2\bZ) \times (\bZ/6\bZ)$ & 4 & $\bigl((1,0),(1,1),(0,2),(0,3)\bigr)$ \cr
\hline
\nr & 4 & $(\bZ/3\bZ) \times (\bZ/3\bZ)$ & 4 & $\bigl((1,0),(1,0),(1,2),(0,1)\bigr)$ \cr
\hline
\nr & 2 & $(\bZ/2\bZ) \times (\bZ/2\bZ)$ & 5 & $\bigl((1,0),(1,0),(1,0),(1,1),(0,1)\bigr)$ \cr
\hline
\nr & 3 & $(\bZ/2\bZ) \times (\bZ/2\bZ)$ & 6 & $\bigl((1,0),(1,0),(1,1),(1,1),(0,1),(0,1)\bigr)$ \cr
\hline
\end{tabular}
\end{center}
\end{table}
\end{question}

Let us give some further details concerning Example~(24). In this case the covers $C_t \to \bP^1$ are branched over $4$~points $t_1,\ldots,t_4$. There are $6$ points above~$t_1$, each with ramification index $e=2$. Likewise, there are $2$~points above~$t_2$, both with $e=6$, there are $4$~points above~$t_3$, each with $e=3$, and finally there are $6$~points above~$t_4$, each with $e=2$. The Hurwitz formula gives $\chi = -24 + 6\cdot 1 + 2\cdot 5 + 4\cdot 2 + 6 \cdot 1 = 6$, so $g=4$. The group ring $Q[A]$, with $A = (\bZ/2\bZ) \times (\bZ/6\bZ)$, is isomorphic to $\bQ^4 \times \bQ(\zeta_3)^4$. Accordingly, the Jacobians~$J_t$ split, up to isogeny, as a product of eight factors. The simple factors of $\bQ[A]$ correspond with the $\Gal(\ol{\bQ}/\bQ)$-orbits in $\Hom(A,\ol{\bQ}^*)$, which in this example may be identified with the set of complex characters $\rho \in \Hom(A,\bC^*)$ taken modulo complex conjugation. Given a complex character $\rho \colon A \to \bC^*$, the corresponding factor $J_{\rho,t}$ of~$J_t$ can be described as follows. We consider the cover $\pi_{\rho,t}\colon C_{\rho,t} \to \bP^1$ with group $\rho(A)$, branched only above the points~$t_i$, with local monodromy $\rho(a_i)$ about~$t_i$. Note that $\rho(a_i)$ may be the identity element of~$\rho(A)$, in which case $t_i$ is not a branch point of~$\pi_{\rho,t}$. Also note that $\pi_{\rho,t}$ is a cyclic cover, which brings us back to the situation considered in Section~\ref{SsT}. Then $J_{\rho,t}$ is the new part of the Jacobian of~$C_{\rho,t}$. With this description, it is easy to verify that:
\begin{itemize}
\item there are five pairs $(\rho,\ol{\rho})$ for which $J_{\rho,t} = 0$; this includes the four real characters, for which $\rho=\ol{\rho}$,
\item there are two pairs $(\rho,\ol{\rho})$ for which $J_{\rho,t}$ is an elliptic curve with CM by an order in~$\bQ(\zeta_3)$,
\item there is one pair $(\rho,\ol{\rho})$ for which $J_{\rho,t}$ is $2$-dimensional, carrying an action by an order in~$\bQ(\zeta_3)$; varying~$t$ these give a family that is isogenous to the family of Jacobians of Example~(5) in Table~\ref{TableExamples}.
\end{itemize}
It follows from this description that the special subvariety~$S(m)$ that contains $Z(m,N,a)$ is $1$-dimensional, and therefore $Z(m,N,a) = S(m)$ is special.

It should be noted that Example~(24) is a sub-family of the family given in Example~(14). To see this, let $D_t$ be the quotient of~$C_t$ modulo the action of $\{1\} \times (\bZ/6\bZ)$, and factor $\pi_t \colon C_t \to \bP^1$ as
\[
C_t \xrightarrow{\; q_t\;} D_t \xrightarrow{\; r_t\;} \bP^1\, .
\]
We have $D_t \cong \bP^1$. The cover~$q_t$ has group $(\bZ/6\bZ)$ and is branched above five points, namely the unique point~$\tilde{t}_2$ of~$D_t$ above~$t_2$, the two points $\tilde{t}_{3,1}$ and~$\tilde{t}_{3,2}$ above~$t_3$, and the two points $\tilde{t}_{4,1}$ and~$\tilde{t}_{4,2}$ above~$t_4$. The local monodromy about these points is given by the $5$-tuple $(2,2,2,3,3)$ in $(\bZ/6\bZ)^5$, which gives exactly the data of Example~(14). The sub-family considered here is given by the constraint that the five branch points on $D_t \cong \bP^1$ do not move freely but form three orbits under the action of $\bZ/2\bZ$ on~$D_t$.

\mosubsection{The Serre-Tate formal group structure and linearity properties}
\label{cancoord}

Let $m \geq 3$. On $\cA_{g,[m]}(\bC)$ we have a metric, obtained from the uniformization by the Siegel space~$\Siegel_g$; see~\ref{AgShimVar}. In \cite{Moonen-LinI} it is proven that an irreducible algebraic subvariety $Z \subset \cA_{g,[m]}$ is a special subvariety if and only if $Z$ is totally geodesic and contains at least one special point. This may be viewed as a characterization of special subvarieties in terms of \emph{linearity properties}.

This characterization has a nice arithmetic analogue, by which it was in fact inspired. Recall that if $k$ is a perfect field of characteristic $p>0$ and $x \in \cA_g(k)$ corresponds to a ppav $(A,\lambda)$ over~$k$ such that $A$ is \emph{ordinary}, the formal completion $\mathfrak{A}_x$ of~$\cA_{g,[m]}$ at the point~$x$ has a canonical structure of a formal torus over the ring of Witt vectors~$W(k)$; see \cite{Katz} or \cite{Messing}, Chap.~5. Using this we can again give meaning to the notion of ``linearity''. Let us elucidate this. Consider a closed irreducible algebraic subvariety $Z \subset \cA_{g,[m],F}$, where $F$ is a number field. Let $\mathcal{Z}$ denote the Zariski closure of~$Z$ inside $\cA_{g,[m]}$ over~$O_F$. If some ordinary point~$x$ as above lies in~$\mathcal{Z}(k)$, the formal completion of~$\mathcal{Z}$ at~$x$ gives a formal subscheme $\mathfrak{Z}_x \subset \mathfrak{A}_x$. It was shown by Noot in~\cite{Noot-Models} that if $Z$ is a special subvariety, the components of the formal subscheme $\mathfrak{Z}_x \subset \mathfrak{A}_x$ over $W(\overline{k})$ are translates of formal subtori of~$\mathfrak{A}_x$ over torsion points. At the cost of excluding finitely many primes of~$O_F$ this may be sharpened to the conclusion that the components of the formal subschemes $\mathfrak{Z}_x$ are formal subtori; see also \cite{Moonen-LinII}, Theorem~4.2(ii).

The converse of Noot's result was proven in~\cite{Moonen-LinII}; the result is that if $Z \subset \cA_{g,[m],F}$ is a closed irreducible algebraic subvariety such that for some ordinary point $x \in \mathcal{Z}(k)$ some component of $\mathfrak{Z}_x \subset \mathfrak{A}_x$ is a translate of a formal subtorus of~$\mathfrak{A}_x$ over a torsion point, $Z$ is a special subvariety. These results again give a characterization of special subvarieties in terms of linearity properties. It was shown in~\cite{Moonen-LinI} that the result over~$\bC$ may reformulated in terms of formal group structures, in a way that makes the analogy between the two situations even clearer.

These results lead us to investigate the structure of the Torelli locus $\Torelli_g \subset \sfA_g$ (or its analogue with a level structure) locally near an ordinary point~$x$ in characteristic~$p$. The identity section of the formal torus~$\mathfrak{A}_x$ gives a lifting of the ppav $(A,\lambda)$ over~$k$ to a polarized abelian scheme $(A^\can,\lambda^\can)$ over the ring of Witt vectors~$W(k)$. This lifting is called the \emph{Serre-Tate canonical lifting} of~$(A,\lambda)$. Note, however, that if $(A,\lambda)$ is the Jacobian of a curve, the canonical lifting  $(A^\can,\lambda^\can)$ over $W(k)$ need not be a Jacobian; see \cite{Dwork.Ogus}, \cite{FO.Sek}. This leads to another view on Expectation~\ref{ZinTorelli}. Indeed, suppose we have a special subvariety $Z \subset \Torelli_g$ with $Z \cap \Torelli_g^\open \neq \emptyset$. Then $Z$ is defined over some number field~$F$, and, with notation as above, we can find ordinary points $x \in \mathcal{Z}(k)$, for finite fields~$k$, such that the corresponding ppav $(A,\lambda)$ is the Jacobian of a curve~$C/k$, and such that the formal completion $\mathfrak{Z}_x \subset \mathfrak{A}_x$ is a union of formal subtori (or even just a formal subtorus). In this situation, the canonical lifting $(A^\can,\lambda^\can)$ is again the Jacobian of a curve, which is a highly nontrivial fact. Indeed, the main results of Dwork and Ogus in~\cite{Dwork.Ogus} are based on the observation that in general, already the first-order canonical lifting, over the ring~$W_2(k)$ of Witt vectors of length~$2$, is no longer a Jacobian.

For the next two questions, let $m\geq 3$ be an integer, and fix a genus~$g$ that is large enough, at least $g>7$. 

\begin{question}
\label{FormalSubtoriQ}
Let $\cT_{g,[m],\bZ} \subset \cA_{g,[m],\bZ}$ denote the scheme-theoretic image of the Torelli morphism $\cM_{g,[m]} \to \cA_{g,[m]}$ over~$\Spec(\bZ)$. Let $k$ be a perfect field of characteristic $p>0$, and suppose $x \in \cT_{g,[m],\bZ}(k)$ is an ordinary point, i.e., the corresponding abelian variety is ordinary. \emph{Is it true that the formal completion $\mathfrak{T}_x \subset \mathfrak{A}_x$ of $\cT_{g,[m],\bZ}$ at the point~$x$ does not contain a formal subscheme $\mathfrak{Z}$, flat over~$W(k)$, of positive dimension, such that $\mathfrak{Z}$ is a formal subtorus of~$\mathfrak{A}_x$?}
\end{question}

A positive answer to this question would confirm Expectation~\ref{ZinTorelli}. It should be noted, however, that what we ask here is stronger than what we need for~\ref{ZinTorelli}. The difference is that in~\ref{FormalSubtoriQ} we do not require the formal subscheme~$\mathfrak{Z}$ to be algebraic, i.e., to be the formal completion of an algebraic subvariety of~$\cA_{g,[m]}$ passing through the point~$x$. The interesting point, however, is that Question~\ref{FormalSubtoriQ} only depends on the formal completion of the Torelli locus at a single ordinary point. One might expect that, in terms of the ``linear structure'' provided by the Serre-Tate structure of a formal torus on~$\mathfrak{A}_x$, the Torelli locus $\mathfrak{T}_x \subset \mathfrak{A}_x$ should be highly non-linear. 

Over $\bC$ we have a question that is similar in spirit.

\begin{question}
\emph{Does the Torelli locus $\cT_{g,[m],\bC} \subset \cA_{g,[m],\bC}$ contain any totally geodesic subvarieties of positive dimension?}
\end{question}

Again this question is stronger than what is needed for Expectation~\ref{ZinTorelli}, as we do not require the totally geodesic subvariety to be algebraic. (In addition we should ask that the subvariety contains at least one CM point.)

Let us mention the paper~\cite{Moller-ST}, in which M\"oller investigates algebraic curves in~$\cM_g$ over~$\bC$ that are totally geodesic with respect to the Teichm\"uller metric (these are called Teichm\"uller curves), such that the image of this curve in~$\cA_g$ is a special subvariety. It is shown in~\cite{Moller-ST} that for $g=2$ and $g\geq 6$ there are no such curves. For $g=3$ and $g=4$ there is precisely one example, corresponding to Examples (7) and~(12) in Table~\ref{TableExamples}.

\mosubsection{An analogy}
\label{analogy}
We can view the boundary of $\cA_g$ in a compactification as the locus of degenerating abelian varieties, but one can also view, working over~$\bZ_p$, say, the space $\cA_g \otimes \bF_p$ as a boundary of $\cA_g \otimes \bZ_p$. In this last case the abelian variety does not degenerate, but the $p$-structure does change. These two points of view have striking similarities. Often geometric questions are settled by studying properties at the boundary, and then to lift back to the interior of the moduli space considered.

Here we want to draw attention to the analogy between the results in~\cite{Dwork.Ogus} and the results of \cite{Fresn.vdP-Unif} and~\cite{Andreatta-CO}. In the first case Dwork and Ogus study ordinary Jacobians in positive characteristic and their Serre-Tate canonical lifts to characteristic zero. In analogy with this, Fresnel and van der Put~\cite{Fresn.vdP-Unif} and Andreatta~\cite{Andreatta-CO} study liftings of Jacobians of degenerate curves. In both cases the authors show that, for $g>3$, these liftings in general do not lie in the Torelli locus. In particular, Andreatta explains in~\cite{Andreatta-CO} that special subvarieties in~$\cA_g$ have a linear structure at the boundary in a toroidal compactification of~$\cA_g$, and he shows that the closure of the Torelli locus is not linear at the boundary.

%%
%%
%% Bibliography
%%

\end{document}

%% file: TorelliMacros
\newcommand{\cA}{\mathcal{A}}

\newcommand{\cG}{\mathcal{G}}
\newcommand{\cH}{\mathcal{H}}

\newcommand{\cM}{\mathcal{M}}

\newcommand{\cP}{\mathcal{P}}
\newcommand{\cT}{\mathcal{T}}

\newcommand{\bA}{\mathbb{A}}
\newcommand{\bC}{\mathbb{C}}
\newcommand{\bF}{\mathbb{F}}
\newcommand{\bG}{\mathbb{G}}

\newcommand{\bP}{\mathbb{P}}
\newcommand{\bQ}{\mathbb{Q}}
\newcommand{\bR}{\mathbb{R}}
\newcommand{\bS}{\mathbb{S}}
\newcommand{\bZ}{\mathbb{Z}}

\newcommand{\gS}{\mathfrak{S}}

\newcommand{\sfA}{\mathsf{A}}

\newcommand{\sfM}{\mathsf{M}}

\newcommand{\ad}{\mathrm{ad}}
\newcommand{\Aut}{\mathrm{Aut}}
\newcommand{\can}{\mathrm{can}}
\newcommand{\CSp}{\mathrm{CSp}}
\newcommand{\ct}{\mathrm{ct}}
\newcommand{\dec}{\mathrm{dec}}
\newcommand{\der}{\mathrm{der}}
\newcommand{\Div}{\mathrm{div}}
\newcommand{\End}{\mathrm{End}}
\newcommand{\eq}{\mathrm{eq}}
\newcommand{\Fil}{\mathrm{Fil}}
\newcommand{\finite}{\mathrm{fin}}
\newcommand{\Gal}{\mathrm{Gal}}
\newcommand{\gen}{\mathrm{gen}}
\newcommand{\GL}{\mathrm{GL}}
\newcommand{\he}{\mathrm{he}}
\newcommand{\Hom}{\mathrm{Hom}}
\newcommand{\id}{\mathrm{id}}
\newcommand{\Inn}{\mathrm{Inn}}

\newcommand{\MT}{\mathrm{MT}}
\newcommand{\mult}{\mathrm{m}}
\newcommand{\new}{\mathrm{new}}
\newcommand{\Nm}{\mathrm{Nm}}
\newcommand{\open}{\circ}

\newcommand{\PGL}{\mathrm{PGL}}
\newcommand{\Pic}{\mathrm{Pic}}
\newcommand{\PSp}{\mathrm{PSp}}
\newcommand{\PSU}{\mathrm{PSU}}
\newcommand{\Res}{\mathrm{Res}}
\newcommand{\SO}{\mathrm{SO}}
\newcommand{\Sp}{\mathrm{Sp}}
\newcommand{\Spec}{\mathrm{Spec}}
\newcommand{\SpecialTor}{\mathcal{ST}}
\newcommand{\SU}{\mathrm{SU}}
\newcommand{\sym}{\mathrm{sym}}

\newcommand{\UU}{\mathrm{U}}

\newcommand{\isomarrow}{\xrightarrow{\sim}}

\newcommand{\Afin}{\bA_\finite}
\newcommand{\Gm}{\bG_\mult}
\newcommand{\Zhat}{\hat{\bZ}}

\newcommand{\Shim}{\mathrm{Sh}}  % general notation for Shimura vars
\newcommand{\Siegel}{\mathfrak{H}}

\newcommand{\Torelli}{\mathsf{T}}
\newcommand{\HElocus}{\mathsf{H}}

\newcommand{\QHS}{\bQ\text{-}\mathsf{HS}}
\let\ol=\overline

\newcount\caseno\caseno=0
\def\nr{\global\advance\caseno by 1({\number\caseno})}

\newcommand{\fracpart}[1]{\left\langle #1\right\rangle}